\documentclass[12pt,review,3p,authoryear]{elsarticle}
\usepackage[psamsfonts]{amssymb}
\usepackage{pict2e}
\usepackage{amsmath,amsthm}
\usepackage{bm}
\usepackage{mathtools}
\usepackage{color}
\usepackage{enumerate}
\usepackage{appendix}
\usepackage{nicefrac}
\usepackage{bbm}
\usepackage{enumitem}
\usepackage{graphicx}
\usepackage{adjustbox}
\usepackage[shortcuts]{extdash} % to hyphenate words with hyphens

\usepackage{relsize}

\usepackage[normalem]{ulem}

\usepackage{tabularx}
\usepackage{multirow}
\usepackage{makecell}
\usepackage{ragged2e}

\usepackage{hyperref}
\usepackage{url}
\usepackage{xr}
\usepackage{cancel}
\usepackage{hhline}
\newcommand{\norm}[1]{\ensuremath{\left|\left|#1\right|\right|}}

\newcommand{\RNum}[1]{\uppercase\expandafter{\romannumeral #1\relax}}

\theoremstyle{definition}
\newtheorem{notation}{Notation}[section]

\theoremstyle{plain}
\newtheorem{theorem}{Theorem}
\newtheorem{lemma}{Lemma}

\theoremstyle{definition}

\newtheorem{comment}{Comment}[section]

\newtheoremstyle{note}% hnamei
{3pt}% hSpace abovei
{3pt}% hSpace belowi
{}% hBody fonti
{}% hIndent amounti
{\itshape}% hTheorem head fonti
{:}% hPunctuation after theorem headi
{.5em}% hSpace after theorem headi
{}% hTheorem head spec (can be left empty, meaning ‘normal’)i
\theoremstyle{note}

\DeclarePairedDelimiter{\abs}{|}{|}
\DeclarePairedDelimiter{\brac}{[}{]}

\DeclarePairedDelimiter{\braces}{\{}{\}}
\DeclarePairedDelimiter{\paren}{\lparen}{\rparen}

\newcommand{\OO}[1]{{\mathcal{O}}\paren*{#1}} 

\newcommand{\MM}[1]{{\Omega}\paren*{#1}}
\newcommand{\other}[1]{\underline{#1}}
\newcommand{\otherr}[1]{\ddot{#1}}
\newcommand{\nzix}[1]{\cancel{0}\paren{#1}}
\newcommand{\nzs}[1]{\abs{\nzix{#1}}}
\newcommand{\orth}[1]{#1^{\perp}}

\newcommand{\prtn}{\mb{P}}
\newcommand{\pbad}{{\tilde{P}}}
\newcommand{\pgood}{{\bar{P}}}
\newcommand{\pz}{{\check{P}}}

\newcommand{\mc}[1]{\ensuremath \mathcal{#1}}
\newcommand{\mb}[1]{\ensuremath \mathbb{#1}}
\newcommand{\real}{\mb{R}}
\newcommand{\extR}{\overline{\real}}

\newcommand{\iden}{\mc{I}_\dimm}
\newcommand{\orthant}{\mb{R}^\dimm_+}

\newcommand{\conj}[1]{\overline{#1}}

\newcommand{\NSpace}[1]{\text{Null}\paren*{#1}}
\newcommand{\notNSpace}[1]{\complement \text{Null}\paren*{#1}}

\newcommand{\cmp}{\mb{C}}
\newcommand{\CmpZero}{\cmp\setminus\braces*{0}}

\newcommand{\Prod}{\mathlarger{\prod}}

\newcommand{\SN}{\text{SN}}
\newcommand{\rMSN}{\text{MSN}}
\newcommand{\uMSN}{\text{CFUSN}}

\newcommand{\TN}{\text{TN}}
\newcommand{\normal}{\text{N}}

\newcommand{\dimm}{K}

\newcommand{\cone}{f_1}
\newcommand{\czero}{f_0}

\newcommand{\czerozz}{\otherr{f}_0}
\newcommand{\Gamz}{\Gamma_0}
\newcommand{\Gamzz}{\otherr{\Gamma}_0}
\newcommand{\Lamz}{\Lambda_0}
\newcommand{\Lamzz}{\otherr{\Lambda}_0}
\newcommand{\Delz}{\Delta_0}
\newcommand{\Delzz}{\otherr{\Delta}_0}
\newcommand{\thez}{\theta_0}
\newcommand{\thezz}{\otherr{\theta}_0}
\newcommand{\muz}{\mu_0}
\newcommand{\muzz}{\otherr{\mu}_0}
\newcommand{\mpuu}{a}

\newcommand{\cf}{CF}
\newcommand{\mgf}{MGF}
\newcommand{\cff}{cf}

\newcommand{\img}{\iota}
\newcommand{\lopi}{L'H\^opital's rule}

\newlist{statement}{enumerate}{1}     % this creates a dedicated counter named 'subtaski'
\setlist[statement,1]{label=\arabic*} % set form of enumeration label here
\title{Identifiability of two-component skew normal mixtures with one known component}

\begin{document}

\author[add1]{Shantanu Jain}
\ead{shajain@indiana.edu}
\author[add2]{Michael Levine}
\ead{mlevins@purdue.edu}
\author[add1]{Predrag Radivojac}
\ead{predrag@indiana.edu}
\author[add3]{Michael W. Trosset}
\ead{mtrosset@indiana.edu}
\address[add1]{Department of Computer Science\\
       Indiana University, Bloomington, Indiana, U.S.A.}
\address[add2]{Department of Statistics\\
       Purdue University, West Lafayette, Indiana, U.S.A.}
\address[add3]{Department of Statistics\\
       Indiana University, Bloomington, Indiana, U.S.A.}
\begin{abstract}
We give sufficient identifiability conditions for estimating mixing proportions in two\-/component mixtures of skew normal distributions with one known component. We consider the univariate case as well as two multivariate extensions: a multivariate skew normal distribution (\rMSN) by \cite{azzalini1996multivariate} and the canonical fundamental skew normal distribution (\uMSN) by \cite{arellano2005fundamental}. The characteristic function of the CFUSN distribution is additionally derived.
\end{abstract}
\maketitle

\section{Introduction}
\label{sec:intro}
\noindent We study identifiability of the mixing proportion for a mixture of two skew normal distributions when one of the components is known. This problem has direct implications for the estimation of mixing proportions given a sample from the mixture and a sample from one of the components. A sample from the mixture is typically collected for a set of objects under study, whereas the component sample is collected for a set of objects verified to satisfy some property of interest. This setting is common in domains where an absence of the property cannot be easily verified due to practical or systemic constraints, e.g., in social networks, molecular biology, etc. In social networks, for example, users may only be allowed to click `like' for a particular product and, thus, the data can be collected only for one of the component samples (a sample from the users who clicked `like') and the mixture (a sample from all users). Accurate estimation of mixing proportions in this setting has fundamental implications for false discovery rate estimation \citep{Storey2002,Storey2003a,Storey2003} and, in the context of classification, for estimating posterior distributions \citep{Ward2009,Jain2016,Jain2016b} and recovering true classifier performance \citep{Menon2015,Jain2017}.

Identifiability and estimation of mixing proportions have been extensively studied \citep{Yakowitz1968, dempster1977maximum, Tallis1982, McLachlan2000}. More recently, the case with one known component has been considered in the nonparametric setting \citep{bordes2006semiparametric, Ward2009, Blanchard2010, Jain2016, Patra2016}. Though the nonparametric formulation is highly flexible, it can also be problematic due to the curse-of-dimensionality issues or when the `irreducibility' assumption is violated \citep{Blanchard2010, Jain2016, Patra2016}. In addition, it is often reasonable in practice to require unimodality of density components, which is difficult to ensure in a nonparametric formulation. To guarantee unimodality of components and allow for skewness, we model the components with a skew normal (SN) family, a generalization of the Gaussian family with good theoretical properties and tractability of inference \citep{genton2004skew}. Although the SN family has been introduced only recently, e.g., see \cite{azzalini1985class, azzalini1986}, and \cite{azzalini1996multivariate}, it has gained practical importance in econometrics and financial domains \citep{genton2004skew}.

Until recently, the literature on identifiability of parametric mixture models emphasized identifiability with respect to a subset of parameters, e.g., cases in which only a single location parameter, or location and scale parameters, can change. Furthermore, most previous results only address the case of univariate mixture distributions. Few studies have considered identifiability of mixtures of general multivariate densities with respect to all of their parameters  \citep{holzmann2006identifiability, browne2015mixture}. 

Our work concerns identifiability with respect to mixing proportions in mixtures of two skew normal distributions with one known component. We  show in Section \ref{sec:iden} that, in this setting, identifiability with respect to mixing proportions is equivalent to  identifiability with respect to all parameters.
We consider both univariate and {\it multivariate}\/ families of skew-normal distributions, establishing identifiability with respect to all of their parameters. We begin with a univariate skew normal family (\SN) introduced by \cite{azzalini1985class}, then extend our results to two forms of multivariate skew normal families (\rMSN\ and \uMSN) introduced by \cite{azzalini1996multivariate} and \cite{arellano2005fundamental}, respectively. These families are further discussed in Section \ref{sec:sn}. Our main contribution is 
%\autoref{thm:SNIden}, \autoref{thm:SNIIden}, and \autoref{thm:CFUSNIden}
Theorems \ref{thm:SNIden}--\ref{thm:CFUSNIden},
which state sufficient conditions for identifiability of the mixing proportion of the mixture with \SN, \rMSN, and \uMSN\ components, respectively. We also derive a concise formula for the characteristic function of \uMSN\ in \ref{sec:CFUSNCF}.

\section{Problem Statement}
\label{sec:probform}
%We are interested in the problem of estimating the mixing proportion of a two-component mixture using a sample from one of the components and a sample from the mixture. %Because of issues of non-identifiability, parameter estimation does not always makes sense \citep{basu1983identifiability}.
%Though this paper primarily deals with mixtures that are univariate, in this section, we develop identifiability theory for the more general multivariate case. 
\noindent Let $\mc{P}_0$ and $\mc{P}_1$ be families of probability density functions (pdfs) on $\mb{R}^{\dimm}$. Let $\mc{F}(\mc{P}_0,\mc{P}_1)$ be a family of pdfs having the form
\begin{align} \label{eq:mixMod}
	f=\alpha f_1 + (1-\alpha)f_0,
\end{align}
%obtained by allowing $f_0$ to vary in $\mc{P}_0$, $f_1$ to vary in $\mc{P}_1$ and $\alpha$ to vary in $(0,1)$. 
where $f_0 \in \mc{P}_0$, $f_1 \in \mc{P}_1$, and $\alpha \in (0,1)$. Densities $f_1$ and $f_0$ will be referred to as component pdfs, $f$ will be referred to as the mixture pdf, and $\alpha$ will be referred to as the mixing proportion. %By definition, 
$\mc{F}(\mc{P}_0,\mc{P}_1)$, therefore, is a family of two-component mixtures. \par

%{\bf This sentence is not needed We are interested in the problem of estimating the mixing proportion $\alpha$ of a two-component mixture $f$ given a sample from one of the components ($f_1$) and a sample from the mixture itself.} 
In this setting, we study identifiability of the density mixture (\ref{eq:mixMod}) with respect to the parameter $\alpha$ %, develop an EM-based algorithm to estimate $\alpha$, and empirically characterize performance of the estimation algorithm. 
when $\mc{P}_0$ and $\mc{P}_1$ are, first, univariate skew-normal and then two different multivariate skew-normal families. All of these distribution families are defined in \citep{genton2004skew}. To do so, we start first with studying some general identifiability conditions in Section \ref{sec:iden}. \par

%Because of identifiability issues, however, parameter estimation does not always makes sense \citep{basu1983identifiability}. This work therefore restricts $\mc{P}_0$ and $\mc{P}_1$ to the families of skew normal distributions

%Let $f_1$ be the component from which we have a sample. From the mathematical viewpoint, this is equivalent (for an infinite sample size) to assuming that $f_1$ is known. Thus, in this section, we will assume that $f_1$ is a fixed {\it known} pdf. This would be same as restricting $\mc{P}_1$ to a singleton set, i.e., $\mc{P}_1=\{f_1\}.$ With a slight abuse of notation, we denote this family by $\mc{F}(\mc{P}_0,f_1).$  Note that $f$ in \autoref{eq:mixMod} can be treated as a pdf parametrized by $\alpha$ and $f_0$. To reflect this parametrization we rewrite $f$ as a function of $\alpha$ and $f_0$, i.e., $f: (0,1) \times \mc{P}_0 \rightarrow \mc{F}(\mc{P}_0,f_1)$ given by $$f(\alpha,f_0)= \alpha f_1 + (1-\alpha)f_0.$$\par

%In this section we discuss identifiability of the mixing proportion in the context of our problem. Our main contribution in this section is \autoref{thm:trnsfrmRatio} that will be used later in the paper.
\section{Identifiabilty of two-component mixtures with a known component}
\label{sec:iden}

\noindent In this section we discuss identifiability of the mixtures in the context of our problem. We will show that the general notion of identifiability is equivalent to identifiability of the mixing proportion (\autoref{lem:twoIden}). However, our main contribution in this section is \autoref{thm:trnsfrmRatio} that gives a useful technique to prove identifiability, tailored to this setting, and will be applied to skew normal mixtures later in the paper. \autoref{lem:twoIden} and \autoref{lem:iden} are restatements of results in \cite{Jain2016}, in terms of densities instead of measures. 

Consider a mixture distribution $f$ from \autoref{eq:mixMod} and let $f_1$ be the known component distribution. This is equivalent to restricting $\mc{P}_1$ to a singleton set, i.e., $\mc{P}_1=\{f_1\}.$ With a minor abuse of notation, we denote the family of mixtures $\mc{F}(\mc{P}_0,\mc{P}_1)$ as $\mc{F}(\mc{P}_0,f_1)$. Note that $f$ in \autoref{eq:mixMod} can be treated as a pdf parametrized by $\alpha$ and $f_0$. To reflect this parameterization, we rewrite $f$ as a function of $\alpha$ and $f_0$, i.e., $f: (0,1) \times \mc{P}_0 \rightarrow \mc{F}(\mc{P}_0,f_1)$ given by $$f(\alpha,f_0)= \alpha f_1 + (1-\alpha)f_0.$$
%A family of distributions $\mc{G}=\braces*{g_\theta:\theta \in \Theta}$ is said to be identifiable if the mapping from $\theta$ to $g_\theta$ is a one-to-one.\footnote{Technically, we require bijection but ignore the obvious ``onto'' requirement for simplicity.} In other words, there exists a unique $\theta$ for every member of $\mc{G}$. Thus, the family of pdfs $\mc{F}(\mc{P}_0,f_1)$ is identifiable if the function $f$ is one-to-one, i.e., $\forall a,b \in (0,1)$ and $\forall h_0,g_0 \in \mc{P}_0,$
A family of distributions $\mc{G}=\braces*{g_\theta:\theta \in \Theta}$ is said to be identifiable if the mapping from $\theta$ to $g_\theta$ is one-to-one.\footnote{Technically, we require bijection but ignore the obvious ``onto'' requirement for simplicity.} Therefore,  $\mc{F}(\mc{P}_0,f_1)$ is identifiable if $\forall a,b \in (0,1)$ and $\forall h_0,g_0 \in \mc{P}_0,$
\begin{align}\label{eq:idenDef}
f(a,h_0) = f(b,g_0) \Rightarrow (a,h_0)=(b,g_0).
\end{align}
The lack of identifiability means that even if $(a,h_0)$ and $(b,g_0)$ are different, the target density $f$ contains no information to tell them apart. If we are only interested in estimating $\alpha$, we need  $\mc{F}(\mc{P}_0,f_1)$ to be identifiable in $\alpha$. That is, $\forall a,b \in (0,1)$ and $\forall h_0,g_0 \in \mc{P}_0,$
\begin{align}\label{eq:idenDefIn}
f(a,h_0) = f(b,g_0) \Rightarrow a=b
\end{align}
Identifiability of $\mc{F}(\mc{P}_0,f_1)$ in $\alpha$ might seem to be a weaker requirement as compared to identifiability of $\mc{F}(\mc{P}_0,f_1)$ in $(\alpha,f_0)$. However, \autoref{lem:twoIden} shows that the two notions of identifiability are equivalent. 
\begin{lemma}
\label{lem:twoIden}
$\mc{F}(\mc{P}_0,f_1)$ is identifiable if and only if $\mc{F}(\mc{P}_0,f_1)$ is identifiable in $\alpha$
\begin{proof}
By definition, identifiability of $\mc{F}(\mc{P}_0,f_1)$ in $\alpha$ is a necessary condition for $\mc{F}(\mc{P}_0,f_1)$ to be identifiable. Now, we prove that it is a sufficient condition as well. Let us assume that $\mc{F}(\mc{P}_0,f_1)$ is identifiable in $\alpha$. Also, suppose that  $\exists a,b \in (0,1)$ and $\exists h_0,g_0 \in \mc{P}_0$ such that $f(a,h_0) = f(b,g_0)$. Then, from the definition of identifiability in $\alpha$, it follows that $a=b$. Therefore, we have $af_1 +(1-a)h_0 = bf_1 + (1-b)g_0$, which implies that $h_0=g_0$. Thus, $\mc{F}(\mc{P}_0,f_1)$ is identifiable.
\end{proof}
\end{lemma}
Consider now the largest possible $\mc{P}_0$, i.e., $\mc{P}_0$ that contains all pdfs in $\mb{R}^\dimm$, except $f_1$ (or any pdf equal to $f_1$ almost everywhere). Then, $\mc{F}(\mc{P}_0,f_1)$ contains all non trivial two component mixtures on $\mb{R}^\dimm$ with $f_1$ as one of the components. \autoref{lem:nonIden} (Section \ref{sec:aux}) shows that this family is not identifiable. Next, we establish the necessary and sufficient condition for identifiability of $\mc{F}(\mc{P}_0,f_1)$.
\begin{lemma}
		\label{lem:iden}
			$\mc{F}(\mc{P}_0,f_1)$ is identifiable if and only if $\mc{F}(\mc{P}_0,f_1) \cap \mc{P}_0=\emptyset.$
\begin{proof}
        First, we will prove that 
        \begin{align}
        \mc{F}(\mc{P}_0,f_1) \cap \mc{P}_0=\emptyset\  \Rightarrow \mc{F}(\mc{P}_0,f_1)\  is\  identifiable. \label{stm:suff}
        \end{align}
        We give a proof by contradiction. Suppose $\mc{F}(\mc{P}_0,f_1) \cap \mc{P}_0=\emptyset$ but that $\mc{F}(\mc{P}_0,f_1)$ is not identifiable. By \autoref{lem:twoIden} $\mc{F}(\mc{P}_0,f_1)$ is not identifiable in $\alpha$. Thus, $\exists a, b \in (0,1)$ and $\exists g_0 ,h_0 \in \mathcal{P}_0$ such that  $f(a,g_0)=f(b,h_0)$, but $a\neq b$. Now, without loss of generality we can assume $a>b$. Therefore, from the equality of $f(a,g_0)$ and $f(b,h_0)$ we obtain, using simple algebra, that $h_0 = \nicefrac{(a-b)}{(1-b)}f_1 + \paren*{1-\nicefrac{(a-b)}{(1-b)}}g_0$. This means, in turn, that $h_0 = f\paren*{\nicefrac{(a-b)}{(1-b)},g_0}$. Because $\nicefrac{(a-b)}{(1-b)} \in (0,1)$ and $g_0 \in \mc{P}_0$, it follows that $h_0 \in \mc{F}(\mc{P}_0,f_1)$. Since $h_0$ has been selected from $\mc{P}_0$ we conclude that $\mc{F}(\mc{P}_0,f_1)\cap \mc{P}_0$ contains $h_0$ and is not empty. This completes the proof of statement \ref{stm:suff}.\\
		Now, we will prove that 
		 \begin{align}
        \mc{F}(\mc{P}_0,f_1)\  is\  identifiable \ \Rightarrow \mc{F}(\mc{P}_0,f_1) \cap \mc{P}_0=\emptyset\  \label{stm:necc}
        \end{align}
        We give a proof by contradiction. Suppose that $\mc{F}(\mc{P}_0,f_1)$ is identifiable but $\mc{F}(\mc{P}_0,f_1) \cap \mc{P}_0 \neq \emptyset$. Let $g_0$ be a common  member of $\mc{P}_0$ and $\mc{F}(\mc{P}_0,f_1)$. As $g_0 \in \mc{F}(\mc{P}_0,f_1)$, it follows that $\exists h_0 \in \mc{P}_0$ and $\exists c \in (0,1)$ such that $g_0 =f(c,h_0)$. Let $a \in (c,1)$ and $b=\nicefrac{(a-c)}{(1-c)}$. As $a,b \in (0,1)$ and $h_0,g_0 \in \mc{P}_0$, it follows that both $f(a,h_0)$ and $f(b,g_0)$ belong in $\mc{F}(\mc{P}_0,f_1)$. We will show that $f(a,h_0)= f(b,g_0).$ Indeed, 
        $f(b,g_0)= bf_1 + (1-b)g_0= bf_1 + (1-b)f(c,h_0)=bf_1 + (1-b)(cf_1+(1-c)h_0)$. This immediately implies that $f(b,g_{0})=f(a,h_{0})$, where $a=b(1-c)+c\in (b,1)$ and is, therefore, greater than $b$. Thus, $(a,h_0) \neq (b,g_0)$. It follows that $\mc{F}(\mc{P}_0,f_1)$ is not identifiable. The lemma follows from statements \ref{stm:suff} and \ref{stm:necc}.     	
\end{proof}
\end{lemma}

The next lemma gives a sufficient condition for identifiability that is mathematically convenient. It relies on the notion of span of a set of functions $\mc{P}$, denoted by $Span(\mc{P})$, that contains all finite linear combinations of functions in $\mc{P}$. That is, $$Span(\mc{P})= \braces*{\sum_{i=1}^{k} a_if_i: k \in \mb{N}, a_i \in \mb{R}, f_i \in \mc{P}}.$$
\begin{lemma}
	\label{thm:trnsfrmRatio}
	Consider the family of pdfs $\mc{F}(\mc{P}_0,f_1)$. Assume that for any pair of pdfs $f_0,g_0 \in \mc{P}_0$ there exists a linear transformation $\Psi$, possibly depending on the choice of $(f_0, g_0)$, that maps any function $f\in \text{Span}(\braces*{f_0,g_0,f_1})$ to a real- or complex-valued function on some domain $S$. We denote $\Psi_{f}$ the value of transformation $\Psi$ of the function $f\in\text{Span}(\braces*{f_0,g_0,f_1})$; thus, $\Psi_{f}$ is a function. We denote $\Psi_{f}(t)$ the value of this function for any $t\in S$. Then, if there exists 
	%such that (\textcolor{red}{PR: below you will say $\Psi_{f_1}(s)$, %so what is the input to this function? Why not define $\Psi$ %formally?})
 a sequence $\braces*{t_n}$ in $S_{f_1} = \braces*{ s \in S: \Psi_{f_1}(s) \neq 0}$ such that
	   \begin{equation}
	             \lim_{n\rightarrow \infty} \frac{\Psi_{f_0}(t_n)}{\Psi_{f_1}(t_n)} = 0 \ \ \text{and} \   \lim_{n\rightarrow \infty} \frac{\Psi_{g_0}(t_n)}{\Psi_{f_1}(t_n)} \notin (-\infty,0), \label{eq:cond}
        \end{equation}
	   %\centerline{or}
        %\item \label{item:f0} there exists sequence $\braces*{t_n}$ in $S_{f_0} = \braces*{ s \in S: \Psi_{f_0}(s) \neq 0}$ such that
        %\begin{align*}
	    %    \lim_{n\rightarrow \infty} \frac{\Psi_{f_1}(t_n)}{\Psi_{f_0}(t_n)} = 0 \ \ \text{and} \   \lim_{n\rightarrow \infty} \frac{\Psi_{g_0}(t_n)}{\Psi_{f_0}(t_n)} \notin (0,1),
       % \end{align*}
    %\end{enumerate}
	$ \mc{F}\paren*{\mc{P}_0, f_1}$ is identifiable.	 
	%\begin{theorem}
%	     The choice of $\Psi$ and $T$ can be different for different pairs of pdfs.
%	\end{theorem}
\begin{proof}
	We give a proof by contradiction. Suppose conditions of the theorem are satisfied but $\mc{F}\paren*{\mc{P}_0, f_1}$ is not identifiable. From \autoref{lem:iden}, it follows that $\mc{F}(\mc{P}_0,f_1) \cap \mc{P}_0 \neq \emptyset$, i.e., there exists a common element in $\mc{F}\paren*{\mc{P}_0, f_1}$ and $\mc{P}_0$, say $f_0$. Because $f_0$ is in $\mc{F}\paren*{\mc{P}_0, f_1}$, there exists $g_0 \in \mc{P}_0$ such that $f_0 = f(a, g_0)$ for some $a \in (0,1)$. Since $f_0$ and $g_0$ are in $\mc{P}_0$, there exists a linear transform $\Psi$ and a sequence $\braces*{t_n}$ satisfying condition (\ref{eq:cond}). It follows that $f_0 = f(a, g_0)= af_1+(1-a)g_0$ and so $ \Psi_{f_0}(t)= a \Psi_{f_1}(t) + (1-a)\Psi_{g_0}(t)$. Now, for all $t\in S_{f_{1}}$ we have $\frac{\Psi_{f_0}(t)}{\Psi_{f_1}(t)} = a + (1-a)\frac{\Psi_{g_0}(t)}{\Psi_{f_1}(t)}$ and consequently, $\lim_{n \rightarrow \infty} \frac{\Psi_{g_0}(t_n)}{\Psi_{f_1}(t_n)} = -\frac{a}{1-a} \in (-\infty,0)$ (contradiction) because $\braces*{t_n}$ satisfies $\lim_{n\rightarrow \infty} \frac{\Psi_{f_0}(t_n)}{\Psi_{f_1}(t_n)}=0$ from condition (\ref{eq:cond}). 
	\end{proof}
\end{lemma}
%For a proof refer \autoref{app1-thm:trnsfrmRatio} in \autoref{app1-app:A}. 
%Later in this paper, we invoke this theorem with two linear transforms, namely the \mgf transform and the \cf\ transform, both of which have $S=\real^\dimm$. \par
We will invoke this lemma later in this paper with two linear transforms; namely, the moment generating function (MGF) transform and the characteristic function (CF) transform. We observe that $S=\real^\dimm$ for both transforms. The main ideas in this lemma (linear transforms and limits) come from Theorem 2 in \cite{teicher1963identifiability} on identifiability of finite mixtures.

\section{Two-component skew normal mixtures}
\label{sec:sn}
%\noindent \autoref{lem:nonIden} shows that the family $\mc{F}(\mc{P}_0,f_1)$, where $\mc{P}_0$ contains all pdfs on $\mb{R}^{\dimm}$\textemdash the family containing all two component mixtures with $f_1$ as one of the components\textemdash is not identifiable. 
%\noindent \autoref{lem:nonIden} shows that the family $\mc{F}(\mc{P}_0,f_1)$ is not identifiable when $\mc{P}_0$ contains all pdfs on $\mb{R}^{\dimm}$. 

\noindent When $\mc{P}_0$ contains all pdfs on $\mb{R}^{\dimm}$, with or without $f_1$, the family $\mc{F}(\mc{P}_0,f_1)$ is not identifiable (\autoref{lem:nonIden}, Section \ref{sec:aux}). It is therefore desirable to choose a smaller family that makes the mixture model identifiable and that is rich enough to model real life data. In this paper, we take a parametric approach. The normal family %\textemdash an obvious choice\textemdash may not be the best option 
presents a limited option since normal mixtures typically require a large number of components to capture asymmetry in real life data. The skew normal family\textemdash an asymmetric family\textemdash provides a convenient alternative in both univatiate and multivariate settings. Thus, we restrict our attention to the two-component mixture families where both unknown and known components are skew normal. Our contribution in this Section are \autoref{thm:SNIden}, \autoref{thm:SNIIden}, and \autoref{thm:CFUSNIden}, that give a rather large identifiable family of two-component skew normal mixtures. A similar approach has been reported by \cite{Ghosal2011} for mixtures of normal and skew normal distributions. Our result, however, results in a much more extensive family $\mc{F}(\mc{P}_0,f_1)$. 
%; in particular, note that the resulting family is much more extensive than the one described in \citet{ghosal2011identifiability}. For convenience, we now provide a more detailed description of the univariate skew normal family as well as its two most common multivariate generalizations. 
Before giving these results, we first introduce the univariate skew normal family as well as its two most common multivariate generalizations. 

\textbf{Univariate skew normal family:} \cite{azzalini1985class} introduced the skew normal (SN) family of distributions as a generalization of the normal family that allows for skewness. It has a location ($\mu$), a scale ($\omega$), and a shape ($\lambda$) parameter, where $\lambda$ controls for skewness. The distribution is right skewed when $\lambda>0$, left skewed when $\lambda<0$, and reduces to a normal distribution when $\lambda=0$. For $X \sim SN(\mu,\omega,\lambda)$, its pdf is given by
%$$f_X(x)=2\omega^{-1}\phi(\omega^{-1}(x-\mu))\Phi(\omega^{-1}\lambda(x-\mu)),\  x \in \mb{R},$$
$$f_X(x)=\frac{2}{\omega}\phi\paren*{\frac{x-\mu}{\omega}}\Phi\paren*{\frac{\lambda(x-\mu)}{\omega}},\  x \in \real,$$
where $\mu,\lambda \in \real$, $\omega \in \real^+$, $\phi$ and $\Phi$ are the probability density function (pdf) and the cumulative distribution function (cdf) of the standard normal distribution $\mc{N}(0,1)$, respectively. \cite{pewsey2003characteristic,genton2004skew,kim2011characteristic} derived the \cf\  and \mgf\ of the \SN\ family (\autoref{tab:CFMGF}).
%\begin{align*}
% \begin{array}{cc}
%    \bar{f}_X(t)= \exp(t\mu + t^2\omega^2/2) \Phi\paren*{\omega\delta t}, &  \hat{f}_X(t)= \exp(it\mu - t^2\omega^2/2) \paren{1 + i\Im(\omega\delta t)}, \ t \in \mb{R},\\
%\end{array}
%\end{align*}

%\textbf{Multivariate skew normal families (\rMSN\ and \uMSN):}
\textbf{Multivariate skew normal families:} \cite{azzalini1996multivariate} proposed an extension of the skew normal family to the multivariate case. This particular generalization has a very useful property in that its marginals are skew normal as well. More recently, several other families of multivariate skew normal distributions have been proposed, as discussed by \cite{lee2013mixtures}. % connects different versions of (msn) families. 
In this paper we consider an alternate parametrization of Azzalini's multivariate skew normal family, denoted by \rMSN. For $X\sim \rMSN(\mu, \Omega,\Lambda)$, pdf of $X$ is   
$$f_X(x|\mu,\Omega,\Lambda)=2\phi_{\dimm}(x-\mu|\Omega)\Phi\paren*{\Lambda'\Omega^{\nicefrac{-1}{2}}(x-\mu)},\ x \in \real^\dimm,$$
where $\Omega$ is a $\dimm \times \dimm$ covariance matrix,  $\mu \in \real^\dimm$ is the location parameter, $\Lambda\in \real^\dimm$ is the shape/skewness parameter, $\phi_\dimm(\cdot|\Omega)$ is the pdf of a $\dimm$-dimensional normal distribution with mean zero and covariance $\Omega$, and $\Phi$ is the cdf of a standard univariate normal. \cite{azzalini1996multivariate,kim2011characteristic} derived the \mgf\ and \cf\ of this distribution (\autoref{tab:CFMGF}).

\citet{lin2009maximum} studied maximum likelihood estimation of finite multivariate skew normal mixtures with another family, the so-called canonical fundamental skew normal distribution (CFUSN), introduced by \cite{arellano2005fundamental}. For $X\sim \uMSN(\mu, \Omega,\Lambda)$, the pdf of $X$ is    
$$f_X(x|\mu,\Omega,\Lambda)=2^\dimm\phi_\dimm(x-\mu|\Omega)\Phi_\dimm\paren*{\Lambda'\Omega^{-1}(x-\mu)|\Delta},\  x \in \real^\dimm$$
where $\Omega$ is a is a $\dimm \times \dimm$ covariance matrix, $\Delta$ is defined in \autoref{tab:altPar}, $\mu \in \mb{R}^\dimm$ is the location parameter, $\Lambda$ is a $\dimm \times \dimm$ shape/skewness matrix\footnote{\cite{arellano2005fundamental} define a more general form of the \uMSN\ family that allows non-square $\Lambda$ matrices.} and $\Phi_\dimm(\cdot|\Delta)$ is the cdf of a $\dimm$-dimensional multivariate normal distribution with zero mean and covariance $\Delta$. The \cf\ and \mgf\ are given in \autoref{tab:CFMGF}. The \mgf\ was obtained from \cite{lin2009maximum}. To the best of our knowledge, the expression for the \cf\ was not available in the literature; we derived it in \autoref{thm:CFUSNCF} (Appendix) for the purposes of this study.

\subsection{Identifiability}
\noindent When $\mc{P}_0$ is some proper/improper subset of the family of univariate skew normal pdfs and $f_1$ is also a univariate skew normal pdf \textemdash concisely written as $\mc{P}_0 \subseteq \braces*{SN(\mu,\omega,\lambda): \mu,\lambda \in \mb{R}, \omega \in \mb{R}^+}$ and $f_1=SN(\mu_1,\omega_1,\lambda_1)$ \textemdash $\mc{F}(\mc{P}_0,f_1)$ contains only two-component univariate skew normal mixtures. \autoref{thm:SNIden} gives a sufficient condition for such a family to be identifiable. 

\begin{table}
%\centering
\caption{Alternate parametrization: The identifiability results and the algorithms are better formulated in terms of the alternate parameters. The table gives the relationship between the alternate and the canonical parameters as well as some other related quantities. Here $\iden$ is a $\dimm \times \dimm$ identity matrix.}
\vspace{2mm}
\begin{adjustbox}{max width=1\textwidth}
\begin{tabularx}{\textwidth}{|c|c|c|c|}
    \cline{1-4}
     \multirow{2}{*}{\textbf{Family}} & \multicolumn{2}{>{\hsize=2\hsize}c|}{ \textbf{Alternate Parametrization}} & \multirow{2}{*}{\textbf{Related Quantities}} \\\cline{2-3}
      & canonical $\rightarrow$ alternate & alternate $\rightarrow$ canonical &  \\ \hhline{====}
      \SN$(\mu,\omega,\lambda)$ & \makecell{ $\Delta=\omega\delta$\\ $\Gamma= \omega^2(1 -\delta^2)$}
            &\makecell{$\omega = \sqrt{\Gamma + \Delta^2 }$\\ $\lambda = sign(\Delta)\sqrt{\nicefrac{\Delta^2}{\Gamma}}$} & $\delta = \frac{\lambda}{\sqrt{1+\lambda^2}}$ \\ \cline{1-4}
            \rMSN$(\mu,\Omega,\Lambda)$& \makecell{$\Delta=\Omega^{\frac{1}{2}}\delta$\\ $\Gamma =\Omega - \Delta\Delta'$}&
            \makecell{$\Omega = \Gamma + \Delta\Delta'$ \\ $\Lambda=\frac{\Omega^{-\frac{1}{2}}\Delta}{\sqrt{1-\Delta'\Omega\Delta}}$}& $\delta= \frac{\Lambda}{\sqrt{1+\Lambda'\Lambda}}$\\ \cline{1-4}
             \uMSN$(\mu,\Omega,\Lambda)$& $\Gamma =\Omega - \Lambda\Lambda'$& $\Omega = \Gamma + \Lambda\Lambda'$ & \makecell{$\Delta =(\iden+\Lambda'\Gamma^{-1}\Lambda)^{-1}$\\
            $= \iden - \Lambda'\Omega^{-1}\Lambda$}\\ \cline{1-4}
\end{tabularx}
\end{adjustbox}
\label{tab:altPar}
\end{table}

\begin{table}
%\centering
\caption{Skew normal families: Expression for the characteristic function and moment generating function. The non-canonical parameters are defined in Table \ref{tab:altPar}. Here $\img$ denotes the imaginary number and $\Im(x)=\int_0^x \sqrt{\nicefrac{2}{\pi}} \exp\paren*{\nicefrac{u^2}{2}}du$. $\Phi$ and $\Phi_\dimm$ denote the cdfs of the standard univariate and multivariate normal distributions, respectively. The $\Lambda_i$ in the expression for \uMSN \ \cf\ is the $i$\textsuperscript{th} column of $\Lambda$. }
\vspace{2mm}
\begin{adjustbox}{center, max width=1\textwidth}
\begin{tabularx}{\textwidth}{|c|c|c|}
    \cline{1-3}
     \textbf{Family} & \textbf{MGF(t)} & \textbf{CF(t)} \\\cline{1-3} \hhline{===}
     \SN$(\mu,\omega,\lambda)$ &$2\exp(t\mu + \nicefrac{t^2\omega^2}{2}) \Phi\paren*{\Delta t}$ & $\exp(\img t\mu - \nicefrac{t^2\omega^2}{2}) \paren{1 + \img\Im(\Delta t)}$ \\ \cline{1-3}
            \rMSN$(\mu,\Omega,\Lambda)$& $2 \exp\paren*{t'\mu + \nicefrac{t' \Omega t}{2}}\Phi(\Delta' t)$ & $\exp\paren*{\img t'\mu - \nicefrac{t' \Omega t}{2}}\paren*{1 + \img\Im( \Delta' t)}$ \\ \cline{1-3}
            \uMSN$(\mu,\Omega,\Lambda)$ & $2^\dimm\exp\braces*{t'\mu + \frac{1}{2}t'\Omega t}\Phi_\dimm\paren*{\Lambda't}$ & $\exp\braces*{\img t'\mu - \frac{1}{2}t'\Omega t}\prod_{i=1}^\dimm\paren*{1+\img\Im(\Lambda_i't)}$\\ \cline{1-3}
\end{tabularx}
\end{adjustbox}
\label{tab:CFMGF}
\end{table}

\begin{notation}[Notation for \autoref{thm:CFUSNIden}, \autoref{lem:samedirection}, \autoref{lem:CFUSNRatio}]
\label{not:CFUSN}
Let $\prtn(U)$ denote a partition defined on a multiset of column vectors $U$ such that column vectors that are in the same direction are in the same set. This relationship is formally defined by the following equivalence relationship 
    \begin{align*}
        t \equiv l \Leftrightarrow ct= l \text{\ for some\ } c\neq 0 \text{\ in\ } \real
    \end{align*}  
   Let $P_C$ denote the canonical vector direction of the vectors in $P\in\prtn(U)$, defined as $P_C= \nicefrac{t}{\norm{t}}$ when a $t \in P$ is not $0$ and $P_C=0$ when a $t\in P$ is $0$. 
   Let $\orth{t}$ be the space orthogonal to a vector $t$.
   Let $\NSpace{M}$ be the null space of matrix $M$. 
   Let $\complement A$ be the complement of a set $A$. 
\end{notation}

\begin{theorem}
	\label{thm:SNIden}
	The family of pdfs $\mc{F}(\mc{P}_0,\cone)$ with $f_1 =\SN(\mu_1,\omega_1,\lambda_1)$ and 
	\begin{align}\label{uniSnfam1}
		 \mathcal{P}_0 = \braces*{SN(\mu,\omega,\lambda) : \Gamma \neq \Gamma_1},
	\end{align}
	($\Gamma$ is defined in Table \ref{tab:altPar}) is identifiable.
	\begin{proof}
			Consider a partition of $\mc{P}_0$ by sets $\mc{P}_0^1$, $\mc{P}_0^2$, defined as follows 
		\begin{align*}
			\mc{P}^1_0 &= \braces*{ \SN(\mu,\omega,\lambda) : \Gamma > \Gamma_1 }\\
			\mc{P}^2_0 &= \braces*{\SN(\mu,\omega,\lambda) : \Gamma < \Gamma_1}\\
		\end{align*}

		We now show that for a given pair of pdfs $\czero,\czerozz$ from $\mc{P}_0$, the conditions of \autoref{thm:trnsfrmRatio} are satisfied. Let $\Gamz,\Delz$ ($\Gamzz,\Delzz$) be the parameters corresponding to $\czero$ ($\czerozz$), as defined in Table \ref{tab:altPar}.
		\begin{itemize}
			\item If $\czero$ is from $\mc{P}^1_0$ ($\Gamz>\Gamma_1$),  we use Lemma \ref{lem:SNRatio} (Statements \ref{item:SNCFRatio0} and \ref{item:SNCFRatio1}) to prove our statement. First, select some $t \neq 0$ in $\mb{R}$. Applying Lemma \ref{lem:SNRatio} (Statements \ref{item:SNCFRatio0} , we obtain
		$\lim_{c \rightarrow \infty} \frac{\cf(ct;\czero)}{\cf(ct;\cone)} = 0$ and $\lim_{c \rightarrow \infty} \frac{\cf(ct;\czerozz)}{\cf(ct;\cone)} \notin (-\infty,0).$ Therefore, the sequence $T=\braces*{t_n}, t_n=nt$ satisfies the conditions of \autoref{thm:trnsfrmRatio}. 
	\item If $\czero$ is from $\mc{P}^2_0$ ($\Gamma_1>\Gamz$),  we use choose Lemma \ref{lem:SNRatio} (Statement \ref{item:SNMGFRatio}) as the basis of our proof. First, we select some $t \neq 0$ in $\mb{R}$ with $\Delz t \leq 0$. Applying Lemma \ref{lem:SNRatio} (Statement \ref{item:SNMGFRatio}), we obtain $\lim_{c \rightarrow \infty} \frac{\mgf(ct;\czero)}{\mgf(ct;\cone)} = 0.$ Moreover, owing to the fact that an mgf is always positive, we know that	$\lim_{c \rightarrow \infty} \frac{\mgf(ct;\czerozz)}{\mgf(ct;\cone)} \notin (-\infty,0).$ The sequence $T=\braces*{t_n}, t_n=nt$ satisfies the conditions of \autoref{thm:trnsfrmRatio}. 
		\end{itemize}
	Thus all the conditions of  \autoref{thm:trnsfrmRatio} are satisfied and consequently $\mc{F}(\mc{P}_0,\cone)$ is identifiable
	\end{proof}
\end{theorem}

\begin{theorem}
	\label{thm:SNIIden}
	The family of pdfs $\mc{F}(\mc{P}_0,\cone)$ with $f_1=\rMSN(\mu_1,\Omega_1,\Lambda_1)$ and  		
	\begin{align*}
		\mathcal{P}_0 = \braces*{\rMSN(\mu,\Omega,\Lambda) : \Gamma \neq \Gamma_1}
	\end{align*}
		($\Gamma$ is defined in Table \ref{tab:altPar}) is identifiable.
		\begin{proof}
	Consider a partition of $\mc{P}_0$ by sets $\mc{P}_0^1$, $\mc{P}_0^2$, defined as follows 
		\begin{align*}
			\mc{P}^1_0 &=  \braces*{\rMSN(\mu,\Omega,\Lambda) : \Gamma \succeq \Gamma_1},\\
			\mc{P}^2_0 &= \mc{P}_0 \setminus \mc{P}_0^1,
		\end{align*}
		where $\succ$ is the standard partial order relationship on the space of matrices. More specifically, $A \succ B$ implies that  $A-B$ is positive definite. Note that $\mc{P}_0^2$ also contains pdfs whose $\Gamma$ matrix is unrelated to $\Gamma_1$ by the partial ordering.\\
		We now show that for a given pair of pdfs $\czero,\czerozz$ from $\mc{P}_0$, the conditions of \autoref{thm:trnsfrmRatio} are satisfied. Let $\Gamz,\Delz$ ($\Gamzz,\Delzz$) be the parameters corresponding to $\czero$ ($\czerozz$), as defined in Table \ref{tab:altPar}.
		\begin{itemize}
			\item If $\czero$ is from $\mc{P}^1_0$,  we choose the characteristic function transform as the linear transform. We pick some $t \in \mb{R}^\dimm$ with $t'(\Gamzz-\Gamma_1)t \neq 0$ and $t'(\Gamz -\Gamma_1)t > 0$; existence of such a $t$ is guaranteed by \autoref{lem:existt}. Applying \autoref{lem:SNIRatio} (Statements \ref{item:SNICFRatio0} and \ref{item:SNICFRatio1}), we obtain $\lim_{c \rightarrow \infty} \frac{\cf(ct;\czero)}{\cf(ct;\cone)}= 0$ and $   	\lim_{c \rightarrow \infty} \frac{\cf(ct;\czerozz)}{\cf(ct;\cone)} \notin (-\infty,0).$ 
		Notice that the sequence $T=\braces*{t_n}, t_n=nt$ satisfies the conditions of \autoref{thm:trnsfrmRatio}. 
		\item If $\czero$ is from $\mc{P}^2_0$,  we choose the moment generating function transform as the linear transform $\Psi$. We pick some $l \neq 0$ in $\mb{R}^k$ such that $l'(\Gamma_1-\Gamz)l > 0$; existence of such an $l$ is guaranteed by $\Gamz \nsucceq \Gamma_1$. If the scalar value $\Delz' l \leq 0$, we choose $t=l$; otherwise, we choose $t=-l$. It is easy to see that $t'(\Gamma_1-\Gamz)t > 0$ and $\Delz' t \leq 0$. Applying \autoref{lem:SNIRatio} (Statement \ref{item:SNIMGFRatio}), we obtain $\lim_{c \rightarrow \infty} \frac{\mgf(ct;\czero)}{\mgf(ct;\cone)} = 0.$ Moreover, owing to the fact that an mgf is always positive, we know that	$\lim_{c \rightarrow \infty} \frac{\mgf(ct;\czerozz)}{\mgf(ct;\cone)} \notin (-\infty,0).$ 
	The sequence $T=\braces*{t_n}, t_n=nt$ satisfies the conditions of \autoref{thm:trnsfrmRatio}. 
		\end{itemize}
	Thus all the conditions of  \autoref{thm:trnsfrmRatio} are satisfied and consequently $\mc{F}(\mc{P}_0,\cone)$ is identifiable
	\end{proof}
\end{theorem}

\begin{theorem}
	\label{thm:CFUSNIden}
	 Let $\theta$ give a concise representation of the \uMSN\ parameters. The family of pdfs $\mc{F}(\mc{P}_0,\cone)$ with $f_1=\uMSN(\theta_1)$ is identifiable 
	when
	%\begin{align*}
	%	\mathcal{P}_0 = \braces*{\uMSN(\theta) : \Gamma \neq \Gamma_1, \uMSN(\theta) \text{\ satisfies condition $\mc{A}$ with\ } \uMSN(\theta_1)},
	%\end{align*}
	\begin{align*}
		\mathcal{P}_0 = \braces*{\uMSN(\theta) : \Gamma \neq \Gamma_1, \Gamma_1-\Gamma \neq  kvv', \text{\ for any\ } v \in \Lambda \text{\ and any\ } k \in \mb{R}^{+}},    
	\end{align*}
	 where $\Gamma$ is defined in Table \ref{tab:altPar}. Here $\Lambda$, in addition to representing the skewness matrix, also represents the multiset containing its column vectors.  
		\begin{proof}
		 First, we define $V(c;\thez,\theta_1,t)=\frac{\exp\paren*{\img c(\muz-\mu_1)'t} \exp\paren*{-\nicefrac{1}{2}c^2t'\paren*{\Gamz-\Gamma_1}t}}{c^{\paren*{\nzs{\Lamz't}-\nzs{\Lambda'_1t}}}},$ where $\nzix{t}=\braces*{i:t[i] \neq 0}$, the set indexes containing non-zero entries of $t$. Note that $V(c;\thez,\theta_1,t)=V(c;\thez,\theta,t)V(c;\theta,\theta_1,t)$, for an arbitrary $\theta$\textemdash a property used multiple times in the proof. We also compute the limit of $V(c;\thez,\theta_1,t)$ as $c \rightarrow \infty$; note that the limit is primarily determined by the sign of the quadratic form $t'(\Gamz-\Gamma_1)t$ and is either $0$ or $\infty$. However, if $t'(\Gamz-\Gamma_1)t=0$, then the limit is determined by the sign of $\nzs{\Lamz't}-\nzs{\Lambda'_1t}$ and is still $0$ or $\infty$; if $\nzs{\Lamz't}-\nzs{\Lambda'_1t}=0$ as well, then $V(c;\thez,\theta_1,t)$ oscillates between $-1$ and $1$ (undefined limit), unless $(\muz-\mu_1)'t=0$, in which case the limit is $1$. We use Notation (\ref{not:CFUSN}) throughout the proof.
		 
		 %Another useful property of $V(c;\thez,\theta_1,t)$ is that
		 %\begin{equation}
		  %      \lim_{c\rightarrow \infty} V(c;\thez,\theta_1,t) \in \CmpZero \Rightarrow  V(c;\thez,\theta_1,t)=1 \label{eq:Vlimit}.
		 %\end{equation}
		 
		 We give a proof by contradiction supposing that the family is not identifiable. Then \autoref{lem:iden} implies that, there exists $\czero$ and $\czerozz$ in  $\mc{P}_0$, such that, with the characteristic function as the linear transform, 
	    \begin{align}
	    \otherr{\cf}_0(ct) & = \mpuu \cf_1(ct) + (1-\mpuu) \cf_0(ct), \quad \forall t \in \real^\dimm,  \forall c \in \real \text{ and }  0<\mpuu<1. \label{eq:CFUSNLin1}
	    \end{align}
	    We will show that \autoref{eq:CFUSNLin1} leads to a contradiction for all possible values of $\czero$ and $\czerozz$ from $\mc{P}_0$.\\
			Consider a partition of $\mc{P}_0$ by sets $\mc{P}_0^1$, $\mc{P}_0^2$, defined as follows 
		\begin{align*}
			\mc{P}^1_0 &=  \braces*{\uMSN(\mu,\Omega,\lambda) : \Gamma_1 \not \succeq \Gamma},\\
			\mc{P}^2_0 &= \mc{P}_0 \setminus \mc{P}_0^1,
		\end{align*}
		where $\succeq$ is the standard partial order relationship on the space of matrices. Precisely, $A \succeq B$ implies that  $A-B$ is positive semi-definite. \\
	Now consider the following cases which cover all the contingencies
	\begin{itemize}
		\item If $\czero$ is from $\mc{P}^1_0$ ($\Gamma_1 \nsucceq \Gamz$) we proceed as follows.
			\autoref{eq:CFUSNLin1} implies that for $\cf_0(ct) \neq 0$
			\begin{align*}
			    &\frac{\otherr{\cf}_0(ct)}{\cf_0(ct)}  = \mpuu\frac{\cf_1(ct)}{\cf_0(ct)} + (1-\mpuu) \\
			    &\quad \Rightarrow \frac{\frac{\otherr{\cf}_0(ct)}{\cf_0(ct)}}{V(c;\theta_1,\thez,t)} = \mpuu \frac{\frac{\cf_1(ct)}{\cf_0(ct)}}{V(c;\theta_1,\thez,t)} + \frac{1-\mpuu}{V(c;\theta_1,\thez,t)} \\
			    &\quad \Rightarrow \frac{1}{V(c;\theta_1,\thezz,t)}\underbrace{\frac{\frac{\otherr{\cf}_0(ct)}{\cf_0(ct)}}{V(c;\thezz,\thez,t)}}_{A} = \underbrace{\mpuu \frac{\frac{\cf_1(ct)}{\cf_0(ct)}}{V(c;\theta_1,\thez,t)}}_{B} + \underbrace{\frac{1-\mpuu}{V(c;\theta_1,\thez,t)}}_{C}
			\end{align*}
	If $\lim_{c\rightarrow\infty} V(c;\theta_1,\thez,t) = \infty$, term (C) goes to $0$ as $c\rightarrow \infty$. Since, applying \autoref{lem:CFUSNRatio} (Statement \ref{item:CFUSNCFAbsRatio}), the limit of term (B) as $c\rightarrow \infty$ exists in $\CmpZero$, so does the limit of the entire RHS and consequently the LHS. It follows that, since limit of term (A) as $c\rightarrow \infty$ exists in $\CmpZero$, $\lim_{c\rightarrow\infty} \frac{1}{V(c;\theta_1,\thezz,t)}$ should also exist in $\CmpZero$ (so that the limit of entire LHS can exist in $\CmpZero$). Summarizing, 
	\begin{equation}
	     \lim_{c\rightarrow\infty} V(c;\theta_1,\thez,t) = \infty \Rightarrow               \lim_{c\rightarrow\infty} \frac{1}{V(c;\theta_1,\thezz,t)} \in \CmpZero
	      \label{eq:CFUSNLin3}
	\end{equation}
	     Now we pick some $t \in \real^\dimm$ with $t'(\Gamz-\Gamma_1)t > 0$ and $t'(\Gamzz-\Gamma_1)t \neq 0$; existence of such a $t$ is guaranteed by $\Gamzz \neq \Gamma_1$ and $\Gamma_1 \nsucceq \Gamz$ as shown in \autoref{lem:existt}. Because $t'(\Gamz-\Gamma_1)t > 0$,  $ \lim_{c\rightarrow\infty} V(c;\theta_1,\thez,t) = \infty$ but $\lim_{c\rightarrow\infty} \frac{1}{V(c;\theta_1,\thezz,t)}$ is either $0$ or $\infty$ as $t'(\Gamzz-\Gamma_1)t \neq 0$, which contradicts \autoref{eq:CFUSNLin3}.
	     
	     	\item If $\czero$ is from $\mc{P}^2_0$, we proceed as follows.
		
			\begin{itemize}
				\item If ($\Gamz = \Gamzz$), we use \autoref{eq:CFUSNLin1} to get
			\begin{align}
			    &\frac{\otherr{\cf}_0(ct)}{\cf_0(ct)}  = \mpuu\frac{\cf_1(ct)}{\cf_0(ct)} + (1-\mpuu) \notag \\
			    &\quad \Rightarrow \frac{\frac{\otherr{\cf}_0(ct)}{\cf_0(ct)}}{V(c;\thezz,\thez,t)} = \mpuu \frac{\frac{\cf_1(ct)}{\cf_0(ct)}}{V(c;\thezz,\thez,t)} + \frac{1-\mpuu}{V(c;\thezz,\thez,t)} \notag\\
			    &\quad \Rightarrow \frac{\frac{\otherr{\cf}_0(ct)}{\cf_0(ct)}}{V(c;\thezz,\thez,t)} = \underbrace{\frac{1-\mpuu}{V(c;\thezz,\thez,t)}}_{A} + \underbrace{\frac{1}{V(c;\thezz,\theta_1,t)}\mpuu \underbrace{\frac{\frac{\cf_1(ct)}{\cf_0(ct)}}{V(c;\theta_1,\thez,t)}}_{B}}_{C} \label{eq:conclusion0}
			\end{align}
			   If $\lim_{c\rightarrow\infty} V(c;\thezz,\theta_1,t) = \infty$, term (C) goes to $0$ as $c\rightarrow \infty$, since the limit of term (B) exists in $\CmpZero$ by  \autoref{lem:CFUSNRatio} (Statement \ref{item:CFUSNCFAbsRatio}). Since, applying \autoref{lem:CFUSNRatio} (Statement \ref{item:CFUSNCFAbsRatio}), the limit of RHS as $c\rightarrow \infty$ exists in $\CmpZero$, so does the limit of the entire LHS and consequently term (A); i.e., $\lim_{c\rightarrow \infty}\frac{1}{V(c;\thezz,\thez,t)} \in \CmpZero$. Summarizing, 
	\begin{equation}
	     \lim_{c\rightarrow\infty} V(c;\thezz,\theta_1,t) = \infty \Rightarrow               \lim_{c\rightarrow\infty} \frac{1}{V(c;\thezz,\thez,t)} \in \CmpZero
	      \label{eq:CFUSNLin4}
	\end{equation}
			    and 
	\begin{align}
	     \lim_{c\rightarrow\infty} V(c;\thezz,\theta_1,t) = \infty &\Rightarrow               \lim_{c\rightarrow \infty} \frac{\frac{\otherr{\cf}_0(ct)}{\cf_0(ct)}}{V(c;\thezz,\thez,t)} = \lim_{c\rightarrow \infty} \frac{1-\mpuu}{V(c;\thezz,\thez,t)} \notag \\
	     &\Rightarrow \Xi(\Lamzz,\Lamz,t) = \lim_{c\rightarrow \infty} \frac{1-\mpuu}{V(c;\thezz,\thez,t)}  \label{eq:CFUSNLin5}\\ 
	     &\tag{from \autoref{lem:CFUSNRatio} (Statement \ref{item:CFUSNCFAbsRatio})}
	\end{align}
	          Now,
	            \begin{align*}
	                t\in \notNSpace{\Gamma_1-\Gamz} & \Rightarrow t'(\Gamma_1-\Gamz)t>0\\
	                & \Rightarrow \lim_{c \rightarrow \infty} V(c;\thezz,\theta_1,t)=\infty\\
	                & \Rightarrow \lim_{c\rightarrow\infty} \frac{1}{V(c;\thezz,\thez,t)} \in \CmpZero \tag{from \autoref{eq:CFUSNLin4}}\\
	                & \Rightarrow \nzs{\Lamz't}- \nzs{\Lamzz't}=0 \text{\ and\ } (\muzz-\muz)'t=0,
	                \end{align*}
	                where the last step follows because $V(c;\thezz,\thez,t)=\frac{\exp\paren*{\img c (\muzz-\muz)'t}}{c^{\paren*{\nzs{\Lamzz't}-\nzs{\Lamz't}}}}$ when $\Gamz=\Gamzz$. Consequently,
	            \begin{equation}
	           \begin{aligned}
	                t\in \notNSpace{\Gamma_1-\Gamz} & \Rightarrow V(c;\thezz,\thez,t)=1\\
	                 & \Rightarrow\Xi(\Lamzz,\Lamz,t) = 1-\alpha  \quad \text{(from \autoref{eq:CFUSNLin5})}
	            \end{aligned} \label{eq:conclusion1}
	            \end{equation}
			 Summarizing, $\forall t \in \notNSpace{\Gamma_1-\Gamz}$,
	            \begin{itemize}
	                \item $\nzs{\Lamz't}- \nzs{\Lamzz't}=0$
	                %\item $(\muzz-\muz)'t=0$
	                \item $\Xi(\Lamzz,\Lamz,t) = 1-\mpuu.$
	            \end{itemize}
	            
	            Since $1-\alpha\neq1$, from \autoref{lem:samedirection} (Statement \ref{item:rneq}), it follows that
	            \begin{itemize}
	              % \item $\forall P \in \prtn(\Lamz \cup \Lamzz), \abs*{P\cap \Lamz} = \abs*{P\cap \Lamzz}$
	               \item $\Gamma_1-\Gamz = kvv'$, for some $v \in \Lamz$ and some $k \in \real^{+}$
	               %\begin{itemize}
	               %\item $\Gamma_1-\Gamz = rP_CP_C'$, for some $r \in \real^{+}$, 
	               %\item $\frac{\prod_{u\in P\cap \Lamzz} \norm{u}}{\prod_{v\in P\cap \Lamz}\norm{v}}=1-\mpuu$,
	               %\item $\frac{\prod_{u\in Q\cap \Lamzz} \norm{u}}{\prod_{v\in Q\cap \Lamz}\norm{v}}=1$, $\forall Q\in \prtn\paren*{\Lamz \cup \Lamzz}\setminus\braces*{P}.$
	               % \end{itemize}               
	            \end{itemize}
	            Thus $\czero \notin \mc{P}_0$ and hence the contradiction.

	       \item If $\Gamz \neq \Gamzz$, 
			
				\autoref{eq:CFUSNLin1} implies that, for $\cf_1(ct) \neq 0$
			\begin{align*}
			    &\frac{\otherr{\cf}_0(ct)}{\cf_1(ct)}  = \mpuu + (1-\mpuu) \frac{\cf_0(ct)}{\cf_1(ct)}\\
			    &\quad \Rightarrow \frac{\frac{\otherr{\cf}_0(ct)}{\cf_1(ct)}}{V(c;\thez,\theta_1,t)}  = \frac{\mpuu}{V(c;\thez,\theta_1,t)} + (1-\mpuu) \frac{\frac{\cf_0(ct)}{\cf_1(ct)}}{V(c;\thez,\theta_1,t)}\\
			  &\quad \Rightarrow \frac{1}{V(c;\thez,\thezz,t)}\underbrace{\frac{\frac{\otherr{\cf}_0(ct)}{\cf_1(ct)}}{V(c;\thezz,\theta_1,t)}}_{A}  =  (1-\mpuu) \underbrace{\frac{\frac{\cf_0(ct)}{\cf_1(ct)}}{V(c;\thez,\theta_1,t)}}_{B} + \underbrace{\frac{\mpuu}{V(c;\thez,\theta_1,t)}}_{C}\\
			\end{align*}
	Notice that if  $\lim_{c\rightarrow\infty} V(c;\thez,\theta_1,t) = \infty$, then term (C) goes to $0$. Since, applying \autoref{lem:CFUSNRatio} (Statement \ref{item:CFUSNCFAbsRatio}), the limit of term (B) as $c\rightarrow \infty$ exists in $\CmpZero$, so is the limit of the entire RHS and consequently the LHS. It follows that, since limit of term (A) as $c\rightarrow \infty$ exists in $\CmpZero$ , $\lim_{c\rightarrow\infty} \frac{1}{V(c;\thez,\thezz,t)}$ should also exist in $\CmpZero$ (so that the limit of entire LHS exists in $\CmpZero$). Summarizing, 
	\begin{equation}
	     \lim_{c\rightarrow\infty} V(c;\thez,\theta_1,t) = \infty \Rightarrow               \lim_{c\rightarrow\infty} \frac{1}{V(c;\thez,\thezz,t)} \in \CmpZero
	      \label{eq:CFUSNLin2}
	\end{equation}
	Now, 
			then we pick some $t \in \mb{R}^\dimm$ with $t'(\Gamma_1-\Gamz)t > 0$ and $t'(\Gamz-\Gamzz)t \neq 0$; existence of such an $t$ is guaranteed by \autoref{lem:existt}. $t'(\Gamma_1-\Gamz)t > 0$ ensures that  $\lim_{c\rightarrow\infty} V(c;\thez,\theta_1,t) = \infty$, but $\lim_{c\rightarrow\infty} V(c;\thez,\thezz,t)$ is either $0$ or $\infty$ as $t'(\Gamz-\Gamzz)t \neq 0$, which contradicts \autoref{eq:CFUSNLin2}.
  \end{itemize}
\end{itemize}
	\end{proof}
\end{theorem}

\begin{comment}[Extension of \autoref{thm:CFUSNIden}]
We speculate that \autoref{thm:CFUSNIden} can be further strengthened by removing the condition $\Gamma_1-\Gamma \neq  kvv'$  (for any $v \in \Lambda$  and any $k \in \mb{R}^{+}$) from the definition of $\mc{P}_0$. Removal of this condition breaks the current proof only in the case when $\Gamma_1 \succeq \Gamz$ and $\Gamz=\Gamzz$. Notice that this case implies that for any $t \in \real^\dimm$ such that $t\in \notNSpace{\Gamma_1-\Gamz}$ satisfies $V(c;\thezz,\theta_1,t) = \MM{c^k\exp\paren*{\nicefrac{1}{2}c^2t'(\Gamma_1-\Gamz)t}}$, for some integer $k$ (from the definition of $V$); $V(c;\thezz,\thez,t)=1$ and $\Xi(\Lamzz,\Lamz,t) = 1-\alpha$ (as shown in \autoref{eq:conclusion1}). These implications reduce \autoref{eq:conclusion0} to
\begin{align*}
  \mpuu \frac{\frac{\cf_1(ct)}{\cf_0(ct)}}{V(c;\theta_1,\thez,t)} &=V(c;\thezz,\theta_1,t)\paren*{\frac{\frac{\otherr{\cf}_0(ct)}{\cf_0(ct)}}{V(c;\thezz,\thez,t)} - \Xi(\Lamzz,\Lamz,t)}\\
  &=\MM{c^k\exp\paren*{\nicefrac{1}{2}c^2t'(\Gamma_1-\Gamz)t}}\Xi(\Lamzz,\Lamz,t)\\
  &\quad \cdot \underbrace{\paren*{\frac{ \prod_{i\in\nzix{\Lamzz't}}1+R_N\paren*{c,(\Lamzz)'_it} -\img \OO{\frac{c}{\exp\paren*{\nicefrac{1}{2}c^2((\Lamzz)'_it)^2}}}}{\prod_{i\in\nzix{\Lamz't}}1+R_N\paren*{c,(\Lamz)'_it} -\img \OO{\frac{c}{\exp\paren*{\nicefrac{1}{2}c^2((\Lamz)'_it)^2}}}} - 1},}_{A} \tag{using \autoref{lem:CFUSNRatio} (Statement \ref{item:CFUSNRatiomm})}
\end{align*}
where, for a positive integer $N$, $R_N(c,x)=\sum_{n=1}^N\frac{(2n-1)!!}{c^{2n}x^{2n}} + \OO{c^{-2(N+1)}}  + \OO{\exp\paren*{-\nicefrac{c^2x^2}{4}}}$ as $c\rightarrow \infty$. Looking at the definition of $R_N$, it seems that the term (A) should be $\MM{c^{k_1}}$ for some negative integer $k_1$, except in a few special cases. This would imply that the RHS is  $\MM{c^{k+k_1}\exp\paren*{\nicefrac{1}{2}c^2t'(\Gamma_1-\Gamz)t}}$ which still goes to $\infty$ as $c\rightarrow \infty$, yet the LHS is in $\CmpZero$, which leads to a contradiction.

\end{comment}

\section{Auxiliary Results}
\label{sec:aux}

\begin{lemma}
\label{lem:nonIden}
If $\mc{P}_0$ contains all pdfs on $\mb{R}^{\dimm}$ except $f_1$, then $\mc{F}=\mc{F}(\mc{P}_0,f_1)$ is not identifiable.
\begin{proof}
Because $\mc{P}_0$ contains all pdfs on $\mb{R}^{\dimm}$ except $f_1$, we have $\mc{F} \subseteq \mc{P}_0$ (note that $f_1 \notin \mc{F}$ either, since $\alpha$ cannot be $1$). Let $a\in(0,1)$ and $b\in(0,a)$, $h_0 \in \mc{P}_0$  and $g_0=f\paren*{\nicefrac{(a-b)}{(1-b)},h_0}$. As $g_0$ is a mixture in $\mc{F}$ and $\mc{F} \subseteq \mc{P}_0$, it follows that $g_0$ is also in $\mc{P}_0$. Consequently, the mixture $f(b, g_0)$ is in $\mc{F}$. Therefore, $f(b,g_0) = bf_1 + (1-b)g_0= bf_1 + (1-b)f\paren*{\nicefrac{(a-b)}{(1-b)},h_0}$; the last expression is equivalent to $f(a,h_{0})$. Thus, we have $f(a,h_0)=f(b, g_0)$. However, $b\neq a$ and hence $\mc{F}$ is not identifiable.
\end{proof}
\end{lemma}

\begin{lemma}
    \label{lem:existt}
    For  $\dimm \times \dimm$ symmetric matrices $A \neq 0$ and $B \neq 0$, if either $A \succeq 0$ or $A \npreceq 0$, then there exists a vector $t \in \real^\dimm$ such that $t'Bt\neq 0$ and $t'At>0$.
    \begin{proof}
        Suppose there does not exist any vector $l \in \real^\dimm$ such that $l'Al>0$. Thus for all  $l \in \real^\dimm$,  $l'Al \leq 0$.  This immediately contradicts $A \npreceq 0$. Hence $A \npreceq 0$ implies that their exists $l \in \real^\dimm$ such that $l'Al>0$. On the other hand, $A \succeq 0$ implies $l'Al \geq 0$ for all   $l \in \real^\dimm$. This, in combination with $l'Al \leq 0$ for all $l \in \real^\dimm$ implies that that $l'Al =0$ for all $l \in \real^\dimm$. This, however, is impossible since $A\neq 0$. Summarizing, there exists $l \in \real^\dimm$ such that $l'Al>0$ when $A\neq 0$ and either of $A \succeq 0$ or $A \npreceq 0$ is true. Now we give a recipe to find $t \in \real^\dimm$ with $t'Bt\neq 0$ and $t'At>0$. Let $l$ be some vector  in $\real^\dimm$ with $l'Al>0$ (existence of $l$ already proved)
        \begin{itemize}
            \item If $l'Bl \neq 0$, then choose $t=l$
            \item else ($l'Bl = 0$) let $l_1 \in \real^\dimm$ be  such that $l_1'Bl_1\neq 0$. Existence of such $l_1$ is guaranteed because $B\neq 0$. We choose $t=l+\epsilon l_1$, where $\epsilon>0$ is picked so that  $t'Bt \neq 0$ and  $t'At>0$. To see that such an $\epsilon$ exists, notice first that $t'Bt= (l+\epsilon l_1)' B (l+\epsilon l_1)=lBl' + 2\epsilon l'_1Bl + \epsilon^2 l'_1  B l_1 =  2\epsilon l'_1Bl + \epsilon^2 l'_1  B l_1 \neq 0$ for any $\epsilon \neq \frac{-2l'_1Bl}{l'_1 B l_1}$. Second, $t'At= (l+\epsilon l_1)' A (l+\epsilon l_1)= l A l' + 2\epsilon l'_1Al + \epsilon^2 l'_1 A l_1>  0$ for a small enough $\epsilon>0$.   Thus picking a small enough $\epsilon \neq  \frac{-2l'_1Bl}{l'_1  B l_1}$ ensures $t'Bt \neq 0$ and  $t'At>0$. 
         \end{itemize}
     \end{proof}
\end{lemma}

\begin{lemma}
\label{lem:samedirection}
Let $U,V$ be $\dimm \times \dimm$ matrices and $S\neq 0$ be a $\dimm \times \dimm$ symmetric positive semi-definite  matrix. Let $U, V$ and $S$ also denote the multiset containing the column vectors of $U,V$ and $S$, respectively and using Notation (\ref{not:CFUSN}), let $\prtn=\prtn(U\cap V)$. Let $\Xi(U,V,t)= \paren*{\img\sqrt{\frac{2}{\pi}}}^{\paren*{ \abs*{\nzix{U't}}-\abs*{\nzix{V't}}}}\frac{\prod_{i\in \nzix{V't}}V'_i t}{\prod_{i\in \nzix{U't}}U'_i t}$, where $\nzix{t}=\braces*{i:t[i] \neq 0}$, $\img$ is the imaginary number and $t\in \real^\dimm$. Assume $\nzs{U't}=\nzs{V't}, \forall t \in \notNSpace{S}$. Then the following statements are true
    \begin{enumerate}
        \item \label{item:sameDirection}   $\abs*{U\cap P}-\abs*{V \cap P} =0, \forall P \in \prtn$; i.e., $P$ has even number of elements with equal contribution from $U$ and $V$.
        \item \label{item:altXi} $\Xi(U,V,l)=\Prod_{P\in \prtn, P_C'l \neq 0}\frac{\prod_{u \in U\cap P} \norm{u}}{\prod_{v\in V\cap P} \norm{v}}, \forall l \in \real^\dimm.$
        \item If $\Xi(U,V,t)=r, \forall t \in \notNSpace{S}$ and some constant $r \in \real$ then 
        \begin{enumerate}
            \item \label{item:rneq} $r\neq 1 \Rightarrow S = kvv'$ for some $v\in V$ and some constant $k>0$.
            \item \label{item:xiNspace} $\Xi(U,V,l)=1, \forall l \in \NSpace{S}$.
        \end{enumerate}
        %\begin{enumerate} \label{item:ratioconst}
        %    \item $\frac{\prod_{u\in P\cap U}u't}{\prod_{v\in P\cap V}v't}=1$,  $\forall t\in \complement \orth{P}_C$, when $s\not \equiv P_C$ for some $s\neq 0$ in $S$ (condition K1),  
        %    \item $\frac{\prod_{u\in P\cap U}u't}{\prod_{v\in P\cap V}v't}=r$, $\forall t\in \complement \orth{P}_C$, when $s \equiv P_C$  for all $s\neq 0$ in $S$ (condition K2). 
        %\end{enumerate}
    \end{enumerate}
\begin{proof}
      First, we partition the elements of $\prtn$ into three sets
     \begin{align*}
         \prtn_0&=\braces*{P\in\prtn: P_C \equiv 0},\\
         \prtn_1&=\braces*{P\in\prtn: P_C \not \equiv s \text{\ for some\ } s \neq 0 \text{\ in\ } S},\\
         \prtn_2&=\prtn\setminus (\prtn_0 \cup \prtn_1).
     \end{align*}
     Notice that $\prtn_0$ is either singleton or empty because all the $0$ vectors in $U\cup V$ are collected in a single component set in $\prtn$. If $\prtn_2\neq \emptyset$, then a vector $w$ in $\pbad \in \prtn_2$ is equivalent to all non-zero column vectors in $S$, which implicitly means that all non-zero column vectors in $S$ are in the same direction (equivalent) and consequently $S$ is rank-$1$ matrix having column vectors (and row vectors  as $S$ is symmetric) equivalent to $w$. In other words, $S$ can be expressed as $S=k_1ww'$ for some constant $k_1>0$ ($k_1>0$ ensures $S$ is positive semi-definite). Summarizing, 
     \begin{equation}
     \prtn_2\neq\emptyset \Rightarrow S=k_1ww', \text{\ for a \ } w \in \pbad \text{\ from $\prtn_2$ and \ } k_1>0 \label{eq:rankone} 
     \end{equation}
     Moreover, any other vector that can appear inside $\prtn_2$ is equivalent to $w$ and consequently $\prtn_2$ is also singleton set (if not empty). These properties are implicitly used in the rest of the proof.      
     
     Next, we show the following result, which will be used multiple times in the proof. \\
     \textbf{(A)}: For a given vector $e$  and a finite multiset of non-zero vectors $M$ in  $\real^\dimm$,  $$e\not\equiv m, \forall m \in M \Rightarrow \exists t \in \real^\dimm \text{\ such that\ } e't=0 \text{\ and\ } m't\neq 0, \forall m \in M.$$
     To prove (A), notice that choosing $t$ from $\orth{e}$ guarantees $e't=0$. Choosing $t$ from $\complement\orth{m}$ ensures $m't \neq 0$. It follows that if the set, $D$, obtained by removing $\orth{m}$, for all $m \in M$, from $\orth{e}$ is non-empty, then any $t \in D$ satisfies both $m't \neq 0$ and $e't=0$. To see that $D$ is indeed non-empty notice that removing  $\orth{m}$'s (finite number of $\dimm-1$ dimensional linear spaces) from  $\orth{e}$ (either $\dimm$ dimensional when $e=0$ or $\dimm-1$ dimensional when $e\neq 0$) reduces it only by Lebesgue measure $0$ set, provided $\orth{e}$ does not coincide with any of the $\orth{m}$'s, guaranteed by $e \not\equiv m$ for all $m \in M$.
     
    Using result (A), we show the existence of two vectors:
    \begin{itemize}
         \item \textbf{$t_0$}: Let $t_0$ be a vector whose existence is shown by using result (A) with $e=0$ and $M=\braces*{P_C: P \in \prtn\setminus\prtn_0}\bigcup\braces*{s}$ for some $s\neq 0$ in $S$. It follows that $P_C't_0\neq 0, \forall P \in \prtn \setminus \prtn_0$ and  $t_0 \in \notNSpace{S}$.
        \item \textbf{$t_{\pgood}$}: For a given $\pgood \in \prtn_1$, let $s\neq 0$ in $S$ be such that $s \not \equiv \pgood_C$ (such an $s$ exists by definition of $\prtn_1$). Let $t_{\pgood}$ be a vector whose existence is shown by using result (A) with $e=\pgood_C$ and $M=\braces*{P_C: P \in \prtn\setminus \paren*{\braces*{\pgood} \cup \prtn_0}}\bigcup\braces*{s}$. It follows that $\pgood_C't_{\pgood}=0$, $P_C't_{\pgood}\neq 0, \forall P \in \prtn \setminus\paren*{\braces*{\pgood} \cup \prtn_0}$ and  $t_\pgood \in \notNSpace{S}$.
    \end{itemize}
     
     To prove the Statement (\ref{item:sameDirection}), we break the argument into three exhaustive cases (picking $P$ from $\prtn_0$ or $\prtn_1$ or $\prtn_2$) as follows
     \begin{enumerate}
         \item \textbf{$\pz \in \prtn_0$}: 
         %From the result (A), there exists a vector $t_0$ such that $P_C't_0=0$, but  $Q_C't_0 \neq 0$ for all $Q \in \prtn \setminus \braces*{P}$ with $Q_C\neq0$ and $s't_0\neq 0$ for some $s\neq 0$ in $S$ 
         %is a non-zero column vector in $S$, ensures that $Q't \neq 0$, for all $Q \in \prtn\setminus\braces*{P}$ and $s't\neq 0$. Existence of $t$ follows from the fact that removing  finite number of $\dimm -1$ dimensional linear spaces ($\orth{Q}$ and $\orth{s}$)  from  $\dimm$ dimensional linear space ($\orth{P}$) gives a non-empty space. 
         Since $t_0 \in \notNSpace{S}$, $\nzs{U't_0}=\nzs{V't_0}$. The only source of $0$'s in $U't_0$ and $V't_0$ are column vectors in $\pz$ and consequently, $\abs*{U\cap \pz}-\abs*{V \cap \pz} =0$ follows.
         \item \textbf{$\pgood \in \prtn_1$}:  
         %Form result (A), there exists a vector $t$ such that $P_C't_P=0$, but  $Q_C't_P \neq 0$ for all $Q \in \prtn \setminus \braces*{P}$ with $Q_C\neq0$ and $s't_P\neq 0$.
         %Picking a vector $t$ from $\orth{P} \setminus \paren*{\bigcup_{Q\in \prtn\setminus\{P\}, Q\neq 0}\orth{Q} \bigcup \orth{s}}$, ensures that $P't=0$, $Q't \neq 0$, for all $Q \in \prtn\setminus\braces*{P}$ such that $Q\neq 0$ and also $s't\neq 0$. Existence of $t$ follows from the fact that the intersection of two distinct $\dimm-1$ spaces has besgue measure $0$ and consequently removing $\orth{Q}$'s and $\orth{s}$ ($\dimm -1$ dimensional linear spaces)  from $\orth{P}$ ($\dimm-1$ dimensional linear space) only reduces the space by a measure $0$ set (none among the $\orth{Q}$'s and $\orth{s}$ is equal to $\orth{P}$). 
          Since $t_{\pgood} \in \notNSpace{S}$, $\nzs{U't_\pgood}=\nzs{V't_\pgood}$.  There are two possibilities for the source of $0$'s in $U't_\pgood$ and $V't_\pgood$: 
         \begin{enumerate}
             \item column vectors in $\pgood$ only (when $\prtn_0=\emptyset$). Thus to satisfy $\nzs{U't_\pgood}=\nzs{V't_\pgood}$, $\abs*{U\cap \pgood}-\abs*{V \cap \pgood} =0$ must be true.
             \item column vectors in $\pgood$ and the only element in $\prtn_0$, $\pz$ (when $\prtn_0$ is singleton). We already know from case (1) that $\abs*{U\cap \pz}-\abs*{V \cap \pz} =0$ is true and consequently, to satisfy $\nzs{U't_\pgood}=\nzs{V't_\pgood}$, $\abs*{U\cap \pgood}-\abs*{V \cap \pgood} =0$ must be true as well.
         \end{enumerate}
          \item \textbf{$\pbad \in \prtn_2$}: Since $\pbad$ is the only element in $\prtn_2$, all other sets $P \in \prtn\setminus\braces*{\pbad}$ belong to either $\prtn_0$ or $\prtn_1$ and are covered by cases (1) and (2); i.e., $\abs*{U\cap P}-\abs*{V \cap P} =0$. As a consequence, $\pbad$, being the only remaining set, $\abs*{U\cap \pbad}-\abs*{V \cap \pbad} =0$ must be true because both $U$ and $V$ have the equal number of column vectors.
         \end{enumerate} 
         This proves Statement (\ref{item:sameDirection}).
     
     To prove Statement  (\ref{item:altXi}), we rewrite the formula for $\Xi(U,V,l), \forall l \in \real^\dimm$ as follows
     \begin{align}
              \Xi(U,V,l)&= \prod_{P\in \prtn\setminus\prtn_0, P_C'l \neq 0}\frac{\prod_{v \in V\cap P} \frac{1}{\img}\sqrt{\frac{\pi}{2}}v'l}{\prod_{u\in U\cap P}\frac{1}{\img}\sqrt{\frac{\pi}{2}}u'l} \notag \\
              &=\prod_{P\in \prtn\setminus\prtn_0, P_C'l \neq 0}\frac{\prod_{v\in V\cap P}\frac{1}{\img}\sqrt{\frac{\pi}{2}}\norm{v}P_C'l}{\prod_{u \in U\cap P} \frac{1}{\img}\sqrt{\frac{\pi}{2}}\norm{u}P_C'l} \notag\\
              &=\prod_{P\in \prtn\setminus\prtn_0, P_C'l \neq 0}\paren*{\img\sqrt{\frac{2}{\pi}} P_C'l}^{\abs*{U\cap P}-\abs*{V \cap P}}\frac{\prod_{v\in V\cap P}\norm{v}}{\prod_{u \in U\cap P} \norm{u}} \notag \\
              &=\prod_{P\in \prtn\setminus\prtn_0, P_C'l \neq 0}\frac{\prod_{v\in V\cap P}\norm{v}}{\prod_{u \in U\cap P} \norm{u}} \tag{because of Statement \ref{item:sameDirection}},
     \end{align}
     which proves Statement (\ref{item:altXi}). 
     Let $t_0$ and $t_\pgood$ (for a given $\pgood \in \prtn_1$) be as defined earlier. Since $t_\pgood$ and $t_0$ are in $\notNSpace{S}$,
     \begin{align}
     \frac{\Xi(U,V,t_0)}{\Xi(U,V,t_\pgood)}=\frac{r}{r} &\Rightarrow \frac{\prod_{P\in \prtn\setminus \prtn_0}\frac{\prod_{v\in V\cap P}\norm{v}}{\prod_{u \in U\cap P} \norm{u}}}{\prod_{P\in \prtn\setminus \paren*{\braces{\pgood} \cup \prtn_0}}\frac{\prod_{v\in V\cap P}\norm{v}}{\prod_{u \in U\cap P} \norm{u}}} = 1.  \notag \\
     &\Rightarrow \frac{\prod_{v\in V\cap \pgood}\norm{v}}{\prod_{u \in U\cap \pgood} \norm{u}}=1 \label{eq:Xione}.
     \end{align}
     Now, 
     \begin{align*}
         \prtn_2 = \emptyset & \Rightarrow \Xi(U,V,l)= \prod_{\pgood\in \prtn_1, \pgood_C'l \neq 0}\frac{\prod_{v\in V\cap \pgood}\norm{v}}{\prod_{u \in U\cap \pgood} \norm{u}}, \forall l \in \real^\dimm\\
            &\Rightarrow \Xi(U,V,l)=1, \forall l \in \real^\dimm \tag{from \autoref{eq:Xione}}\\
             &\Rightarrow \Xi(U,V,t)=1, \forall t \in \notNSpace{S}. 
     \end{align*}
       Thus,
       \begin{align*}
           r\neq 1 &\Rightarrow \prtn_2 \neq \emptyset\\
                   &\Rightarrow S = k_1ww', \text{\ for a\ } w \in \pbad \text{\ from \ } \prtn_2 \text{\ and some \ } k_1>0 \tag{from \autoref{eq:rankone}}\\
                   &\Rightarrow S = kvv', \text{\ for some \ } v \in V \cap \pbad \text{\ and some\ } k>0;
       \end{align*}
               existence of $v$ is justified by Statement (\ref{item:sameDirection}) and the fact that $\pbad$ is non-empty (contains $w$). This proves Statement (\ref{item:rneq}).
        
        To prove Statement (\ref{item:xiNspace}), notice that if $\prtn_2 \neq \emptyset$, then $\forall l \in \NSpace{S}$ and for the only set $\pbad \in \prtn_2$, $\pbad'_C l=0$ (from the definition of $\prtn_2$). It follows that $\forall l \in \NSpace{S}$       
               \begin{align*}
                   \Xi(U,V,l)&= \prod_{\pgood\in \prtn_1, \pgood_C'l \neq 0}\frac{\prod_{v\in V\cap \pgood}\norm{v}}{\prod_{u \in U\cap \pgood} \norm{u}}, \forall l \in \real^\dimm \tag{because either $\prtn_2=\emptyset$ or $\pbad_C'l=0$ for the only $\pbad \in \prtn_2$}\\
            &=1. 
               \end{align*}
\end{proof}
\end{lemma}

\begin{notation}[Landau's notation]
\label{not:landau}
We use Landau's asymptotic notation in the next few lemmas, defined as follows. For real-valued functions $g$ and $h$ defined on some subset of $\real$, $g(c)=\OO{h(c)}$ as $c\rightarrow \infty$ if $\limsup_{c\rightarrow \infty} \abs*{\frac{g(c)}{h(c)}} < \infty$ and $g(c)=\MM{h(c)}$ as $c\rightarrow \infty$ if $\limsup_{c\rightarrow \infty} \abs*{\frac{g(c)}{h(c)}} > 0$. 
\end{notation}

\begin{lemma}
	\label{lem:SNRatio}
	Consider two univariate skew normal distributions, $\SN(\mu,\omega,\lambda)$ and $\SN(\other{\mu},\other{\omega},\other{\lambda})$. Let $\Gamma,\Delta$ be related to $\omega$ and $\lambda$ as given in Table \ref{tab:altPar}. Let $c\in \mb{R}$ and $t \in \mb{R}\setminus\braces*{0}$.
		\begin{enumerate}
		\item Let $\cf$ and $\other{\cf}$ be the characteristic functions corresponding to the two distributions (refer Table \ref{tab:CFMGF}).
		    \begin{enumerate}
		        \item \label{item:SNCFRatio0}
		    	$$ \Gamma - \other{\Gamma} > 0 \Rightarrow \lim_{c \rightarrow \infty} \frac{\cf(ct)}{\other{\cf}(ct)} = 0$$
	        \item \label{item:SNCFRatio1}
		      $$ \Gamma - \other{\Gamma} \neq 0 \Rightarrow \lim_{c \rightarrow \infty} \frac{\cf(ct)}{\other{\cf}(ct)} \in \braces*{-\infty,0,\infty},$$
		        provided the limit exists in $\extR$ (the extended real number line).
		    \end{enumerate}
        \item \label{item:SNMGFRatio} Let $\mgf$ and $\other{\mgf}$ be the moment generating functions corresponding to the two distributions (refer Table \ref{tab:CFMGF}). For $\Delta t \leq 0$,
         $$\other{\Gamma} - \Gamma > 0 \Rightarrow \lim_{c \rightarrow \infty} \frac{\mgf(ct)}{\other{\mgf}(ct)} = 0.$$
		\end{enumerate}
\begin{proof} 
    Here, we use Landau's $\OO{\cdot}$ and $\MM{\cdot}$ notation, defined in Notation (\ref{not:landau}).\\
     \textbf{Statement \ref{item:SNCFRatio0}}: Instead of working directly with $\frac{\cf(ct)}{\other{\cf}(ct)}$, which can be complex, we circumvent the complication by working with the ratio's absolute value squared, which is always real. Multiplying the ratio with its conjugate, we obtain an expression of its absolute value squared as follows:
		\begin{align*}
			\abs*{\frac{\cf(ct)}{\other{\cf}(ct)}}^2&=\frac{\cf(ct)}{\other{\cf}(ct)}\conj{\paren*{\frac{\cf(ct)}{\other{\cf}(ct)}}}\notag\\
			&= \frac{\cf(ct)\conj{\paren*{\cf(ct)}}}{\other{\cf}(ct)\conj{\paren*{\other{\cf}(ct)}}}\tag{property of complex conjugate of a fraction}\\
			&= \frac{\exp\paren*{-c^2\omega^2_0t^2}\paren*{1+\paren*{\Im\paren*{c\Delta t}}^2}}{\exp\paren*{-c^2t^2\omega^2_1t}\paren*{1+\paren*{\Im\paren*{c\other{\Delta}t}}^2}}
		\end{align*}
		Consider the ratio $\frac{1+\paren*{\Im\paren*{c\Delta t}}^2}{1+\paren*{\Im\paren*{c\other{\Delta}t}}^2}$ from the previous expression. Using the asymptotic upper-bound (for the numerator) and lower bound (for the denominator), obtained in \autoref{lem:asymptoticBounds} (Statement \ref{item:ImOO} and \ref{item:ImMM}),  we get
	\begin{align*}
			\frac{1+\paren*{\Im\paren*{c\Delta t}}^2}{1+\paren*{\Im\paren*{c\other{\Delta}t}}^2} &= \OO{c^2\exp\paren*{c^2\paren*{\Delta t}^2 -c^2\paren*{\other{\Delta}t}^2}}\\
		&=\OO{c^2\exp\paren*{c^2 t^2\paren*{\Delta^2- \other{\Delta}^2}}}
	\end{align*}
		
Thus,
\begin{align}
	\abs*{\frac{\cf(ct)}{\other{\cf}(ct)}}^2&= \exp\paren*{-c^2t^2(\omega^2-\other{\omega}^2)}\OO{c^2\exp\paren*{c^2t^2\paren*{\Delta^2 - \other{\Delta}^2}}} \notag\\
    &= \OO{c^2\exp\paren*{- c^2t^2\paren*{\paren*{\omega^2-\Delta^2} - \paren*{\other{\omega}^2 -\other{\Delta}^2}}}}\notag\\
	&=\OO{c^2\exp\paren*{- c^2t^2\paren*{\Gamma - \other{\Gamma} }}} \label{eq:SNCfRatioUb}
	\end{align}
Consequently, 
\begin{align*}
& \lim_{c \rightarrow \infty} 	\abs*{\frac{\cf(ct)}{\other{\cf}(ct)}}^2 =0, \text{\ when\ } \Gamma - \other{\Gamma}>0
\end{align*}
and 
\begin{align*}
&\lim_{c \rightarrow \infty} \frac{\cf(ct)}{\other{\cf}(ct)} =0, \text{\ when\ } \Gamma - \other{\Gamma}>0
\end{align*}
follows.

\textbf{Statement \ref{item:SNCFRatio1}}  Similar to the derivation of the asymptotic upper-bound for the ratio in Equation \ref{eq:SNCfRatioUb}, we derive the asymptotic lower-bound by using Lemma \ref{lem:asymptoticBounds} (Statement \ref{item:ImOO} and\ref{item:ImMM});
\begin{align*}
	\abs*{\frac{\cf(ct)}{\other{\cf}(ct)}}^2	&= \MM{\nicefrac{1}{c^2}\exp\paren*{c^2t^2\paren*{\other{\Gamma} - \Gamma }}}
\end{align*}
 Consequently, 
\begin{align*}
& \lim_{c \rightarrow \infty} 	\abs*{\frac{\cf(ct)}{\other{\cf}(ct)}}^2 =\infty, \text{\ when\ } \other{\Gamma} - \Gamma>0 
\end{align*}
and 
\begin{align*}
&\lim_{c \rightarrow \infty} \frac{\cf(ct)}{\other{\cf}(ct)} \in\braces*{-\infty,\infty}, \text{\ when\ } \other{\Gamma} - \Gamma >0 
\end{align*}
follows, provided the limit exists in $\extR$.
Combining the result with Statement \ref{item:SNICFRatio0} proves Statement \ref{item:SNICFRatio1}.

\textbf{Statement \ref{item:SNMGFRatio}}
From the definition of \SN\  \mgf\ (Table \ref{tab:CFMGF}) we get
		\begin{align*}
			\frac{\mgf(ct)}{\other{\mgf}(ct)}&= \exp\paren*{c(\mu-\other{\mu})t - \frac{c^2}{2}t^2(\other{\omega}^2-\omega^2)}\frac{\Phi\paren*{c\Delta t}}{\Phi\paren*{c\other{\Delta}t}}
		\end{align*}
		Consider the ratio $\frac{\Phi\paren*{c\Delta t}}{\Phi\paren*{c\other{\Delta}t}}$ from the previous expression. We apply the asymptotic upper-bound (for the numerator) and lower bound (for the denominator), obtained in \autoref{lem:asymptoticBounds} (Statement \ref{item:PhiOO} and \ref{item:PhiMM}). Because $\Delta t \leq 0$, the asymptotic upper-bound is applicable.  
		\begin{align*}
		\frac{\Phi\paren*{c\Delta t}}{\Phi\paren*{c\other{\Delta}t}} &= \OO{c\exp\paren*{-\frac{c^2}{2}\paren*{(\Delta t)^2 - (\other{\Delta}t)^2}}}\\
		&=\OO{c\exp\paren*{-\frac{c^2}{2}t^2\paren*{\Delta^2 - \other{\Delta}^2}}} 
	\end{align*}
Thus,
\begin{align*}
	\frac{\mgf(ct)}{\other{\mgf}(ct)}&= \exp\paren*{c(\mu-\other{\mu})t- \frac{c^2}{2}t^2(\other{\omega}^2-\omega^2)}\OO{c\exp\paren*{-\frac{c^2}{2}t^2\paren*{\Delta^2 - \other{\Delta}^2}}} \notag\\
	&= \OO{c\exp\paren*{c(\mu-\other{\mu})t - \frac{c^2}{2}t^2\paren*{\paren*{\other{\omega}^2-\other{\Delta}^2}- \paren*{\omega^2-\Delta^2} }}}\\
	&= \OO{c\exp\paren*{c(\mu-\other{\mu})t - \frac{c^2}{2}t^2\paren*{\other{\Gamma} -\Gamma}}}
	\end{align*}
Because $c^2$ term dominates the $c$ term in the exponential above, the asymptotic upper-bound goes to $0$ when $\other{\Gamma} -\Gamma > 0$, irrespective of the relation between $\mu$ and $\other{\mu}$. Consequently, 
\begin{align*}
\lim_{c \rightarrow \infty} 	\frac{\mgf(ct)}{\other{\mgf}(ct)} = 0, \text{\ when\ } \other{\Gamma} - \Gamma>0.
\end{align*}
	\end{proof}
\end{lemma}

\begin{lemma}
	\label{lem:SNIRatio}
	Consider two $\dimm-$dimensional Skew Normal distributions, $\rMSN(\mu,\Omega,\Lambda_0)$ and $\rMSN(\other{\mu},\other{\Omega},\Lambda_1)$. Let $\Gamma,\Delta$ be related to $\Omega$ and $\Lambda$ as given in Table \ref{tab:altPar}. Let $c\in \mb{R}$ and $t \in \mb{R}^\dimm$.
		\begin{enumerate}
		\item Let $\cf$ and $\other{\cf}$ be the characteristic functions corresponding to the two distributions (refer Table \ref{tab:CFMGF}).
		    \begin{enumerate}
		        \item \label{item:SNICFRatio0}
		    	$$ t'\paren*{\Gamma - \other{\Gamma}}t > 0 \Rightarrow \lim_{c \rightarrow \infty} \frac{\cf(ct)}{\other{\cf}(ct)} = 0$$
	        \item \label{item:SNICFRatio1}
		      $$ t'\paren*{\Gamma - \other{\Gamma}}t \neq 0 \Rightarrow \lim_{c \rightarrow \infty} \frac{\cf(ct)}{\other{\cf}(ct)} \in \braces*{-\infty,0,\infty},$$
		        provided the limit exists in $\extR$ (the extended real number line).
		    \end{enumerate}
        \item \label{item:SNIMGFRatio} Let $\mgf$ and $\other{\mgf}$ be the moment generating functions corresponding to the two distributions (refer Table \ref{tab:CFMGF}). For $\Delta' t \leq 0$,
         $$t'\paren*{\other{\Gamma} - \Gamma}t > 0 \Rightarrow \lim_{c \rightarrow \infty} \frac{\mgf(ct)}{\other{\mgf}(ct)} = 0.$$
		\end{enumerate}
\begin{proof} 
    Here, we use Landau's $\OO{\cdot}$ and $\MM{\cdot}$ notation, defined in Notation (\ref{not:landau}).\\
     \textbf{Statement \ref{item:SNICFRatio0}}: We use the approach in \autoref{lem:SNRatio}. The expression for the squared absolute value of the characteristic function ratio, obtained by multiplying the ratio with its conjugate, is given by
		\begin{align*}
			\abs*{\frac{\cf(ct)}{\other{\cf}(ct)}}^2&= \frac{\exp\paren*{-c^2t'\Omega t}\paren*{1+\paren*{\Im\paren*{c\Delta't}}^2}}{\exp\paren*{-c^2t'\other{\Omega}t}\paren*{1+\paren*{\Im\paren*{c\other{\Delta}'t}}^2}}
		\end{align*}
		Consider the ratio $\frac{1+\paren*{\Im\paren*{c\Delta't}}^2}{1+\paren*{\Im\paren*{c\other{\Delta}'t}}^2}$ from the previous expression. Using the asymptotic upper-bound (for the numerator) and lower bound (for the denominator), obtained in \autoref{lem:asymptoticBounds} (Statement \ref{item:ImOO} and \ref{item:ImMM}),  we get
	\begin{align*}
			\frac{1+\paren*{\Im\paren*{c\Delta't}}^2}{1+\paren*{\Im\paren*{c\other{\Delta}'t}}^2} &= \OO{c^2\exp\paren*{c^2\paren*{\Delta't}^2 - c^2\paren*{\other{\Delta}'t}^2}}\\
		&=\OO{c^2\exp\paren*{c^2 t'\paren*{\Delta\Delta' - \other{\Delta}\other{\Delta}'}t}}, 
	\end{align*}
Thus,
\begin{align}
	\abs*{\frac{\cf(ct)}{\other{\cf}(ct)}}^2&= \exp\paren*{-c^2t'(\Omega-\other{\Omega})t}\OO{c^2\exp\paren*{c^2t'\paren*{\Delta\Delta' - \other{\Delta}\other{\Delta}'}t}} \notag\\
    &= \OO{c^2\exp\paren*{- c^2t'\paren*{\paren*{\Omega-\Delta\Delta'} - \paren*{\other{\Omega} -\other{\Delta}\other{\Delta}'}}t}}\notag\\
	&=\OO{c^2\exp\paren*{- c^2\paren*{t'\paren*{\Gamma - \other{\Gamma}}t}}} \label{eq:SNICfRatioUb}
	\end{align}
Consequently, 
\begin{align*}
& \lim_{c \rightarrow \infty} 	\abs*{\frac{\cf(ct)}{\other{\cf}(ct)}}^2 =0, \text{\ when\ } t'\paren*{\Gamma - \other{\Gamma}}t>0
\end{align*}
and 
\begin{align*}
&\lim_{c \rightarrow \infty} \frac{\cf(ct)}{\other{\cf}(ct)} =0, \text{\ when\ } t'\paren*{\Gamma - \other{\Gamma}}t>0
\end{align*}
follows.

\textbf{Statement \ref{item:SNICFRatio1}} Similar to the derivation of the asymptotic upper-bound for the ratio in Equation \ref{eq:SNICfRatioUb}, we derive the asymptotic lower-bound by using Lemma \ref{lem:asymptoticBounds} (Statement \ref{item:ImOO} and \ref{item:ImMM});
\begin{align*}
	\abs*{\frac{\cf(ct)}{\other{\cf}(ct)}}^2	&= \MM{\frac{1}{c^2}\exp\paren*{c^2\paren*{t'\paren*{\other{\Gamma} - \Gamma}t}}}
\end{align*}
Consequently, 
\begin{align*}
& \lim_{c \rightarrow \infty} 	\abs*{\frac{\cf(ct)}{\other{\cf}(ct)}}^2 =\infty, \text{\ when\ } t'\paren*{\other{\Gamma} - \Gamma}t>0 
\end{align*}
and 
\begin{align*}
&\lim_{c \rightarrow \infty} \frac{\cf(ct)}{\other{\cf}(ct)} \in\braces*{-\infty,\infty}, \text{\ when\ } t'\paren*{\other{\Gamma} - \Gamma}t>0 
\end{align*}
follows, provided the limit exists in $\extR$.
Combining the result with Statement \ref{item:SNICFRatio0} proves Statement \ref{item:SNICFRatio1}.

\textbf{Statement \ref{item:SNIMGFRatio}}
From the definition of \rMSN\  \mgf\ (Table \ref{tab:CFMGF}) we get
		\begin{align*}
			\frac{\mgf(ct)}{\other{\mgf}(ct)}&= \exp\paren*{c(\mu'-\other{\mu}')t - \frac{c^2}{2}t'(\other{\Omega}-\Omega)t}\frac{\Phi\paren*{c\Delta't}}{\Phi\paren*{c\other{\Delta}'t}}
		\end{align*}
		Consider the ratio $\frac{\Phi\paren*{c\Delta't}}{\Phi\paren*{c\other{\Delta}'t}}$ from the previous expression. We apply the asymptotic upper-bound (for the numerator) and lower bound (for the denominator), obtained in \autoref{lem:asymptoticBounds} (Statement \ref{item:Phi}). Because $\Delta' t \leq 0$, the asymptotic upper-bound is applicable.  
		\begin{align*}
		\frac{\Phi\paren*{c\Delta't}}{\Phi\paren*{c\other{\Delta}'t}} &= \OO{c\exp\paren*{-\frac{c^2}{2}\paren*{(\Delta't)^2 - (\other{\Delta}'t)^2}}}\\
		&=\OO{c\exp\paren*{-\frac{c^2}{2}t'\paren*{\Delta\Delta' - \other{\Delta}\other{\Delta}'}t}} 
	\end{align*}
Thus,
\begin{align*}
	\frac{\mgf(ct)}{\other{\mgf}(ct)}&= \exp\paren*{c(\mu'-\other{\mu}')t- \frac{c^2}{2}t'(\other{\Omega}-\Omega)t}\OO{c\exp\paren*{-\frac{c^2}{2}t'\paren*{\Delta\Delta' -\other{\Delta}\other{\Delta}'}t}} \notag\\
	&= \OO{c\exp\paren*{c(\mu'-\other{\mu}')t - \frac{c^2}{2}t'\paren*{\paren*{\other{\Omega}-\other{\Delta}\other{\Delta}'} - \paren*{\Omega -\Delta\Delta'} }t}}\\
	&= \OO{c\exp\paren*{c(\mu'-\other{\mu}')t - \frac{c^2}{2}\paren*{t'\paren*{\other{\Gamma} -\Gamma}t }}}
	\end{align*}
Because $c^2$ term dominates the $c$ term in the exponential above, the asymptotic upper-bound goes to $0$, irrespective of the relation between $\mu$ and $\other{\mu}$. Consequently, 
\begin{align*}
\lim_{c \rightarrow \infty} 	\frac{\mgf(ct)}{\other{\mgf}(ct)} = 0, \text{\ when\ } t'\paren*{\other{\Gamma} - \Gamma}t>0.
\end{align*}
	\end{proof}
\end{lemma}

\begin{lemma}
	\label{lem:CFUSNRatio}
	Consider two $\dimm-$dimensional Skew Normal distributions, $\uMSN(\mu,\Omega,\Lambda)$ and $\uMSN(\other{\mu},\other{\Omega},\other{\Lambda})$. 
	Let $\Gamma,\Delta$ be related to $\Omega$ and $\Lambda$ as given in Table \ref{tab:altPar}. Let $c\in \mb{R}$ and $t \in \mb{R}^\dimm$.
		\begin{enumerate}
		\item Let $\cf$ and $\other{\cf}$ be the characteristic functions corresponding to the two distributions (refer Table \ref{tab:CFMGF}). Let $V(c;\theta,\other{\theta},t)=\frac{\exp\paren*{\img(\mu-\other{\mu})'t} \exp\paren*{-\nicefrac{1}{2}c^2t'\paren*{\Gamma-\other{\Gamma}}t}}{c^{\paren*{\nzs{\Lambda't}-\nzs{\other{\Lambda}'t}}}}$, where $\nzix{t}=\braces*{i:t[i] \neq 0}$. Let $\Xi(U,V,t)= \paren*{\img\sqrt{\frac{2}{\pi}}}^{\paren*{ \abs*{\nzix{U't}}-\abs*{\nzix{V't}}}}\frac{\prod_{i\in \nzix{V't}}V'_i t}{\prod_{i\in \nzix{U't}}U'_i t}$, where $\img$ is the imaginary number and $U_i$ ($V_i$) is the $i$\textsuperscript{th} column of $U$ ($V$). Then
		\begin{enumerate} 
		    \item  \label{item:CFUSNCFAbsRatio}
		     %\begin{align*}
			%\lim_{c\rightarrow\infty} \frac{c^{\paren*{\nzs{\Lambda't}-\nzs{\other{\Lambda}'t}}}\abs*{\frac{\cf(ct)}{\other{\cf}(ct)}}}{\exp\paren*{-\nicefrac{1}{2}c^2t'\paren*{\Gamma-\other{\Gamma}}t}} &= \paren*{\frac{2}{\pi}}^{\nicefrac{1}{2}\paren*{ \nzs{\Lambda't}-\nzs{\other{\Lambda}'t}}}\prod_{i=1}^\dimm\frac{\paren*{\other{\Lambda}'_it}^{ \indc{\other{\Lambda}'_it\neq0}}}{\paren*{\Lambda'_it}^{ \indc{\Lambda'_it\neq0}}},
			%\end{align*}
			\begin{align*}
			\lim_{c\rightarrow\infty} \frac{\frac{\cf(ct)}{\other{\cf}(ct)}}{V(c;\theta,\other{\theta},t)} &=  \Xi(\Lambda,\other{\Lambda},t).
			\end{align*}
			%\item \label{item:CFUSNRatiooo} 
		%	$$ \frac{\frac{\cf(ct)}{\other{\cf}(ct)}}{V(c;\theta,\other{\theta},t)} -  \Xi(\Lambda,\other{\Lambda},t) = \oo{\nicefrac{1}{c}}$$
		   % where $\Xi(U,V,t)= \paren*{\img\sqrt{\frac{2}{\pi}}}^{\paren*{ \abs*{\nzix{U't}}-\abs*{\nzix{V't}}}}\frac{\prod_{i\in \nzix{V't}}V'_i t}{\prod_{i\in \nzix{U't}}U'_i t}$ and $\nzix{t}=\braces*{i:t[i] \neq 0}$ and $U_i$ ($V_i$) is the $i$\textsuperscript{th} column of $U$ ($V$).
		  %When $\Lambda't$ is not a permutation of $\other{\Lambda}'t$, for some positive integer $k$
		   % $$ \frac{\frac{\cf(ct)}{\other{\cf}(ct)}}{V(c;\theta,\other{\theta},t)} -  \Xi(\Lambda,\other{\Lambda},t) = \mm{c^{-k}}$$
		    
		     \item \label{item:CFUSNRatiomm}    Using Landau's $\OO{\cdot}$ notation, defined in Notation (\ref{not:landau}),
		     $$ \frac{\frac{\cf(ct)}{\other{\cf}(ct)}}{V(c;\theta,\other{\theta},t)}-\Xi(\Lambda,\other{\Lambda},t)=\Xi(\Lambda,\other{\Lambda},t) \paren*{\frac{ \prod_{i\in\nzix{\Lambda't}}1+R_N\paren*{c,\Lambda'_it} -\img \OO{\frac{c}{\exp\paren*{\nicefrac{1}{2}c^2(\Lambda'_it)^2}}}}{\prod_{i\in\nzix{\other{\Lambda}_i't}}1+R_N\paren*{c,\other{\Lambda}'_it} -\img \OO{\frac{c}{\exp\paren*{\nicefrac{1}{2}c^2(\other{\Lambda}'_it)^2}}}} - 1},$$
		    where, for a positive integer $N$, $R_N(c,x)=\sum_{n=1}^N\frac{(2n-1)!!}{c^{2n}x^{2n}} + \OO{c^{-2(N+1)}}  + \OO{\exp\paren*{-\nicefrac{c^2x^2}{4}}}$ as $c\rightarrow \infty$.
		\end{enumerate}
		   
		\end{enumerate}
\begin{proof}
     		\begin{align}
			\frac{\frac{\cf(ct)}{\other{\cf}(ct)}}{V(c;\theta,\other{\theta},t)}&= c^{\paren*{\nzs{\Lambda't}-\nzs{\other{\Lambda}'t}}}\frac{\exp\paren*{-\nicefrac{1}{2}c^2t'\paren*{\Omega-\other{\Omega}}t}}{\exp\paren*{-\nicefrac{1}{2}c^2t'\paren*{\Gamma-\other{\Gamma}}t}}\prod_{i=1}^\dimm\frac{1+\img \Im\paren*{c\Lambda'_it}} {1+\img\Im\paren*{c\other{\Lambda}'_it}}\notag\\
			&= \frac{c^{\nzs{\Lambda't}}}{c^{\nzs{\other{\Lambda}'t}}}\frac{\exp\paren*{\nicefrac{1}{2}c^2t'\other{\Lambda}\other{\Lambda}'t}}{\exp\paren*{\nicefrac{1}{2}c^2t'\Lambda\Lambda't}}\prod_{i=1}^\dimm\frac{1+\img \Im\paren*{c\Lambda'_it}}{1+\img \Im\paren*{c\other{\Lambda}'_it}}\notag\\
			&= \frac{c^{\nzs{\Lambda't}}}{c^{\nzs{\other{\Lambda}'t}}}\prod_{i=1}^\dimm\frac{\frac{1+\img\Im\paren*{c\Lambda'_it}}{\exp\paren*{\nicefrac{1}{2}c^2(\Lambda'_it)^2}}}{\frac{1+\img\Im\paren*{c\other{\Lambda}'_it}}{\exp\paren*{\nicefrac{1}{2}c^2(\other{\Lambda}'_it)^2}}}\notag\\
			&=\frac{ \prod_{i\in\nzix{\Lambda't}}\frac{c\paren*{1+\img\Im\paren*{c\Lambda'_it}}}{\exp\paren*{\nicefrac{1}{2}c^2(\Lambda'_it)^2}}}{\prod_{i\in\nzix{\other{\Lambda}'t}}\frac{c\paren*{1+\img\Im\paren*{c\other{\Lambda}'_it}}}{\exp\paren*{\nicefrac{1}{2}c^2(\other{\Lambda}'_it)^2}}}\notag\\
				&=\frac{ \prod_{i\in\nzix{\Lambda't}}\frac{c}{\exp\paren*{\nicefrac{1}{2}c^2(\Lambda'_it)^2}}+\img \frac{c}{\exp\paren*{\nicefrac{1}{2}c^2(\Lambda'_it)^2}}\Im\paren*{c\Lambda'_it}}{\prod_{i\in\nzix{\other{\Lambda}_i't}}\frac{c}{\exp\paren*{\nicefrac{1}{2}c^2(\other{\Lambda}'_it)^2}}+\img \frac{c}{\exp\paren*{\nicefrac{1}{2}c^2(\other{\Lambda}'_it)^2}}\Im\paren*{c\other{\Lambda}'_it}} \label{eq:CFUSNCFRatio}
			\end{align}
       Using \autoref{lem:asymptoticBounds} (Statement \ref{item:Imlim}), we get
       \begin{align*}
			\lim_{c\rightarrow \infty}\frac{\frac{\cf(ct)}{\other{\cf}(ct)}}{V(c;\theta,\other{\theta},t)} &= \frac{\prod_{i\in \nzix{\Lambda't}}\paren*{\img\sqrt{\frac{2}{\pi}}\frac{1}{\Lambda'_it}}}{\prod_{i\in \nzix{\other{\Lambda}'t}}\paren*{\img\sqrt{\frac{2}{\pi}}\frac{1}{\other{\Lambda}'_it}}}\\
			%&= \paren*{\frac{2}{\pi}}^{\paren*{ \nzs{\Lambda't}-\nzs{\other{\Lambda}'t}}}\prod_{i=1}^\dimm\frac{\paren*{\other{\Lambda}'_it}^{2 \indc{\other{\Lambda}'_it\neq0}}}{\paren*{\Lambda'_it}^{2 \indc{\Lambda'_it\neq0}}}.
			&= \Xi(\Lambda,\other{\Lambda},t).
		\end{align*}
			This proves Statement (\ref{item:CFUSNCFAbsRatio}).\\
		Using \autoref{eq:CFUSNCFRatio}, 	
		\begin{align*}
		       \frac{\frac{\cf(ct)}{\other{\cf}(ct)}}{V(c;\theta,\other{\theta},t)}-\Xi(\Lambda,\other{\Lambda},t) &=\frac{ \prod_{i\in\nzix{\Lambda't}}\frac{c}{\exp\paren*{\nicefrac{1}{2}c^2(\Lambda'_it)^2}}+\img \frac{c}{\exp\paren*{\nicefrac{1}{2}c^2(\Lambda'_it)^2}}\Im\paren*{c\Lambda'_it}}{\prod_{i\in\nzix{\other{\Lambda}_i't}}\frac{c}{\exp\paren*{\nicefrac{1}{2}c^2(\other{\Lambda}'_it)^2}}+\img \frac{c}{\exp\paren*{\nicefrac{1}{2}c^2(\other{\Lambda}'_it)^2}}\Im\paren*{c\other{\Lambda}'_it}} - \Xi(\Lambda,\other{\Lambda},t)\\
		       &=\frac{ \prod_{i\in\nzix{\Lambda't}}\frac{c}{\exp\paren*{\nicefrac{1}{2}c^2(\Lambda'_it)^2}}+\img \sqrt{\nicefrac{2}{\pi}}\frac{1}{\Lambda't}\paren*{1+R_N\paren*{c,\Lambda'_it}}}{\prod_{i\in\nzix{\other{\Lambda}_i't}}\frac{c}{\exp\paren*{\nicefrac{1}{2}c^2(\other{\Lambda}'_it)^2}}+\img \sqrt{\nicefrac{2}{\pi}}\frac{1}{\other{\Lambda}'t}\paren*{1+R_N\paren*{c,\other{\Lambda}'_it}}} - \Xi(\Lambda,\other{\Lambda},t) \tag{Using \autoref{lem:asymptoticBounds} Statement(\ref{item:Impoly})}\\
		       &=\Xi(\Lambda,\other{\Lambda},t) \frac{ \prod_{i\in\nzix{\Lambda't}}1+R_N\paren*{c,\Lambda'_it} -\img \OO{\frac{c}{\exp\paren*{\nicefrac{1}{2}c^2(\Lambda'_it)^2}}}}{\prod_{i\in\nzix{\other{\Lambda}_i't}}1+R_N\paren*{c,\other{\Lambda}'_it} -\img \OO{\frac{c}{\exp\paren*{\nicefrac{1}{2}c^2(\other{\Lambda}'_it)^2}}}} - \Xi(\Lambda,\other{\Lambda},t)\\
		       &=\Xi(\Lambda,\other{\Lambda},t) \paren*{\frac{ \prod_{i\in\nzix{\Lambda't}}1+R_N\paren*{c,\Lambda'_it} -\img \OO{\frac{c}{\exp\paren*{\nicefrac{1}{2}c^2(\Lambda'_it)^2}}}}{\prod_{i\in\nzix{\other{\Lambda}_i't}}1+R_N\paren*{c,\other{\Lambda}'_it} -\img \OO{\frac{c}{\exp\paren*{\nicefrac{1}{2}c^2(\other{\Lambda}'_it)^2}}}} - 1}
	    \end{align*}
	This proves Statement (\ref{item:CFUSNRatiomm}).
	\end{proof}

\end{lemma}

%\begin{lemma}
%    \label{lem:asymptoticProperty}
%    Let $a(x),b(x)$ be complex functions such that $\lim_{x\rightarrow\infty}a(x)=a\neq 0$ and $\lim_{x\rightarrow\infty}b(x)=b \neq 0$. Further, let $a(x)-a=\oo{u(x)}$, $a(x)-a=\mm{v(x)}$, $b(x)-b=\oo{u(x)}$ and $b(x)-b=\mm{v(x)}$. Then
%    \begin{enumerate}
%        \item \label{item:ooPropMult} $a(x)b(x)-ab=\oo{u(x)}$.
%        \item \label{item:mmPropMult} $a(x)b(x)-ab=\mm{v(x)}$.
%        \item \label{item:ooPropDiv} $a(x)b(x)-ab=\oo{u(x)}$.
%        \item \label{item:mmPropDiv} $a(x)b(x)-ab=\mm{v(x)}$.
%    \end{enumerate}
%    \begin{proof}
%         \begin{align*}
%             a(x)b(x)-ab &= a(x)b(x)-ab(x)+ab(x)-ab\\
%                       &= b(x)(a(x)-a) + a(b(x)-b)\\
%                       &=b(x)\oo{u(x)} + a\oo{u(x)}\\
%                       &=\OO{1}\oo{u(x)} + a\oo{u(x)}\\
%                       &=\oo{u(x)}
%        \end{align*}
%   \end{proof}
%\end{lemma}

\begin{lemma}
\label{lem:asymptoticBounds}
	Let $\Phi$ be the standard normal cdf and $\Im(x) =  \int_0^x \sqrt{\nicefrac{2}{\pi}} \exp \paren*{\nicefrac{u^2}{2}} du$. Let $x$ be finite. Then, using Landau's $\OO{\cdot}$ and $\MM{\cdot}$ notation, defined in Notation (\ref{not:landau}), as $c\rightarrow \infty$, 
		\begin{enumerate}
		\item \label{item:Phi} 
			\begin{enumerate}
			\item \label{item:PhiMM} For all $x\in\real$
			\begin{align*}
			    \Phi(cx)=\MM{\frac{\exp\paren*{-\nicefrac{c^2x^2}{2}}}{c}} 
			\end{align*}
			\item \label{item:PhiOO} When $x\leq 0$,
			\begin{align*}
			    \Phi(cx)=\OO{exp\paren*{-\nicefrac{c^2x^2}{2}}}  
			\end{align*}
		\end{enumerate}
	    \item \label{item:Im}  %$\Imexpr(c,x)=c\frac{1+\img\Im(cx)}{\exp(\nicefrac{c^2x^2}{2})}$. Then, for $x\neq 0$,
		    \begin{enumerate}
		        \item For all $x\neq 0$, \label{item:Imlim} 
		        \begin{align*}
	                \lim_{c\rightarrow \infty}
	                 \frac{c\Im(cx)}{\exp\paren*{\frac{c^2x^2}{2}}}= \sqrt{\nicefrac{2}{\pi}} \frac{1}{x}
	            \end{align*}
	            \item For all $x\neq 0$, \label{item:Impoly}
	           $$ \frac{c\Im(cx)}{\exp\paren*{\frac{c^2x^2}{2}}}= \sqrt{\frac{2}{\pi}} \frac{1}{x} \brac*{1 + \sum_{n=1}^N\frac{(2n-1)!!}{c^{2n}x^{2n}} + \OO{c^{-2(N+1)}}  + \OO{\exp\paren*{-\nicefrac{c^2x^2}{4}}}}, $$
	           where $!!$ is defined recursively for an integer $a\geq -1$ as follows
	           $$a!!= \begin{cases}
	                    a(a-2)!! & \text{when $a\neq0$ and $a\neq -1$}\\ 
	                    1 & \text{when $a=0$ or $a=-1$}
	                    \end{cases}$$
	            % \item %\label{item:Imexpmm} 
	            %and for an arbitrary $\epsilon>0$ 
	            %$$\Imexpr(c,x)- \img\sqrt{\nicefrac{2}{\pi}} \nicefrac{1}{x} = \mm{\exp(-\epsilon c)}.$$
	            \item For all $x \in \real$, \label{item:ImOO}
	            $$1+ \paren*{\Im(cx)}^2 = \OO{\exp(c^2x^2)}$$
	            \item For all $x\in \real$, \label{item:ImMM}
		        $$1+ \paren*{\Im(cx)}^2 = \MM{\frac{\exp(c^2x^2)}{c}}$$
		    \end{enumerate}
    \end{enumerate}
	\begin{proof}\hspace*{\fill}\\
	\textbf{Statement \ref{item:Phi}}:
	
	Consider the function, 
	\begin{align*}
		g(c) &= c\frac{\Phi(cx)}{\exp\paren*{-\frac{c^2x^2}{2}}},\\
	\end{align*}
	To prove the statements, we first  derive the limits of $g(c)$ as $c \rightarrow \infty$. 
	%When $x\geq0$,
	%\begin{equation}
	%	\begin{aligned} \label{eq:PhiExpRatio1}
	%	    \lim_{c \rightarrow \infty} g(c) &= \begin{cases}
	%	        \frac{1}{2} & \text{ when } x=0, \\
	%	         \infty & \text{ when } x>0. 
	%	    \end{cases}
    %    \end{aligned}
    %\end{equation}
	To evaluate the limit, when $x<0$, we apply the \lopi\ because both the numerator and denominator go to $0$. To this end we take the limit of ratio of the derivative of the numerator and the denominator w.r.t $c$ and applying Leibniz integral rule, we get
	\begin{align}
	     \lim_{c\rightarrow \infty}g(c)&=\lim_{c\rightarrow \infty}\frac{\frac{d}{dc}\Phi(cx)}{\frac{d}{dc}\frac{\exp\paren*{-\frac{c^2x^2}{2}}}{c}} \notag \\
	     &= \lim_{c\rightarrow \infty} \frac{\frac{x}{\sqrt{2\pi}} \exp\paren*{-\frac{c^2x^2}{2}}}{\frac{\paren*{(-cx^2)c-1}\exp\paren*{-\frac{c^2x^2}{2}}}{c^2}} \notag \\
	     &= \lim_{c\rightarrow \infty} \frac{\frac{x}{\sqrt{2\pi}} }{\paren*{-x^2-\nicefrac{1}{c^2}}} \notag \\
	      &= \frac{-1}{\sqrt{2\pi}x} \label{eq:PhiExpRatio2}
	\end{align}
	Thus, for $x<0$, $\Phi(cx)=\OO{\frac{\exp\paren*{-\nicefrac{c^2x^2}{2}}}{c}}$ and consequently, $\Phi(cx)=\OO{exp\paren*{-\nicefrac{c^2x^2}{2}}}$, which also holds true when $x=0$. Thus $\Phi(cx)=\OO{exp\paren*{-\nicefrac{c^2x^2}{2}}}$ when $x\leq0$, which proves Statement (\ref{item:PhiOO}).	Moreover, it follows from \autoref{eq:PhiExpRatio2} that, for $x<0$, $\Phi(cx)=\MM{\frac{\exp\paren*{-\nicefrac{c^2x^2}{2}}}{c}}$ and since  it is true for $x\geq0$ as well (because $\Phi(0)=\nicefrac{1}{2}$ and $\Phi(cx)$ approaches $1$ when $x>0$), $\Phi(cx)=\MM{\frac{\exp\paren*{-\nicefrac{c^2x^2}{2}}}{c}}$ for all $x\in\real$. This proves Statement (\ref{item:PhiMM}).

	%\begin{align*}
	%    \Phi(cx)=\OO{\exp\paren*{-\nicefrac{c^2x^2}{2}}}
	%\end{align*}
    %Summarizing, from \ref{eq:PhiExpRatio1} and \ref{eq:PhiExpRatio2}, 			
    %$$  \lim_{c \rightarrow \infty}	\frac{\Phi(cx)}{\exp\paren*{-\frac{c^2x^2}{2}}} < \infty \  \text{when}\ x \leq 0,$$
	%Thus $\Phi(cx)=\OO{\exp \paren*{-\frac{c^2x^2}{2}}}$ when $x \leq 0$\\
	%Also, from \ref{eq:PhiExpRatio1} and \ref{eq:PhiExpRatio2}, 			
    %$$  \lim_{c \rightarrow \infty}	\frac{\Phi(cx)}{\exp\paren*{-(1+\epsilon)\frac{c^2x^2}{2}}} > 0 \  \text{when}\ \epsilon > 0, \ \text{for all}\ x$$
	%Thus $\Phi(cx)=\MM{\exp \paren*{-(1+\epsilon)\frac{c^2x^2}{2}}}$.\\
    \textbf{Statement \ref{item:Im}}: Performing integration by parts on $\Im(cx)$, for $x\neq 0$ gives
    
     \begin{align*}
        \Im(x)&=\int_0^{x} \sqrt{\frac{2}{\pi}}  \exp\paren*{\nicefrac{u^2}{2}} du\\
            &= \int_0^{\nicefrac{x}{\sqrt{2}}} \sqrt{\frac{2}{\pi}}  \exp\paren*{\nicefrac{u^2}{2}} du +\int_{\nicefrac{x}{\sqrt{2}}}^{x} \sqrt{\frac{2}{\pi}}  \exp\paren*{\nicefrac{u^2}{2}} du\\
            &=\sqrt{\frac{2}{\pi}} \brac*{\int_0^{\nicefrac{x}{\sqrt{2}}}  \exp\paren*{\nicefrac{u^2}{2}} du +  \int_{\nicefrac{x}{\sqrt{2}}}^{x} \frac{1}{u} \frac{d}{du}\paren*{\exp\paren*{\nicefrac{u^2}{2}}}du} \\
            &= \sqrt{\frac{2}{\pi}} \brac*{\int_0^{\nicefrac{x}{\sqrt{2}}}  \exp\paren*{\nicefrac{u^2}{2}} du - \frac{\exp\paren*{\nicefrac{x^2}{4}}}{2^{\nicefrac{-1}{2}}x} + \frac{\exp\paren*{\nicefrac{x^2}{2}}}{x}  + \int_{\nicefrac{x}{\sqrt{2}}}^{x} \frac{\exp\paren*{\nicefrac{u^2}{2}}}{u^2}du}\\
             &=  \sqrt{\frac{2}{\pi}} \brac*{ \int_0^{\nicefrac{x}{\sqrt{2}}}  \exp\paren*{\nicefrac{u^2}{2}} du - \frac{\exp\paren*{\nicefrac{x^2}{4}}}{2^{\nicefrac{-1}{2}}x} + \frac{\exp\paren*{\nicefrac{x^2}{2}}}{x} + \int_{\nicefrac{x}{\sqrt{2}}}^{x} \frac{1}{u^3}\frac{d}{du}\exp\paren*{\nicefrac{u^2}{2}}du}\\
             &= \sqrt{\frac{2}{\pi}} \brac*{\int_0^{\nicefrac{x}{\sqrt{2}}}  \exp\paren*{\nicefrac{u^2}{2}} du - \frac{\exp\paren*{\nicefrac{x^2}{4}}}{2^{\nicefrac{-1}{2}}x}-\frac{\exp\paren*{\nicefrac{x^2}{4}}}{2^{\nicefrac{-3}{2}}x^3} +\frac{\exp\paren*{\nicefrac{x^2}{2}}}{x} + \frac{\exp\paren*{\nicefrac{x^2}{2}}}{x^3} + \int_{\nicefrac{x}{\sqrt{2}}}^{x} 3\frac{\exp\paren*{\nicefrac{u^2}{2}}}{u^4}du}\\
             &= \sqrt{\frac{2}{\pi}} \left[  \frac{\exp\paren*{\nicefrac{x^2}{2}}}{x} +\sum_{n=1}^N\frac{(2n-1)!!  \exp\paren*{\nicefrac{x^2}{2}}} {x^{2n+1}} + (2N+1)!!\int_{\nicefrac{x}{\sqrt{2}}}^x\frac{  \exp\paren*{\nicefrac{u^2}{2}}}{u^{2(N+1)}} du \right.\\
             &\left. -\sum_{n=0}^N\frac{(2n-1)!!\exp\paren*{\nicefrac{x^2}{4}}}{\sqrt{2^{-(2n+1)}}x^{2n+1}} +\int_0^{\nicefrac{x}{\sqrt{2}}}  \exp\paren*{\nicefrac{u^2}{2}} du  \right] \\
              &= \sqrt{\frac{2}{\pi}}\frac{\exp\paren*{\nicefrac{x^2}{2}}}{x} \left[ 1 + \sum_{n=1}^N\frac{(2n-1)!!}{x^{2n}}  + (2N+1)!! x \frac{\int_{\nicefrac{x}{\sqrt{2}}}^x \frac{  \exp\paren*{\nicefrac{u^2}{2}}}{u^{2(N+1)}} du}{\exp\paren*{\nicefrac{x^2}{2}}}  \right.\\
              &\left. -\sum_{n=0}^N\frac{(2n-1)!!\exp\paren*{-\nicefrac{x^2}{4}}}{\sqrt{2^{-(2n+1)}}x^{2n}} + x\frac{\int_0^{\nicefrac{x}{\sqrt{2}}} \exp\paren*{\nicefrac{u^2}{2}} du}{\exp\paren*{\nicefrac{x^2}{2}}} \right]
    \end{align*}
    Thus, 
    \begin{align*}
        \Im(cx) &= \sqrt{\frac{2}{\pi}}\frac{\exp\paren*{\nicefrac{c^2x^2}{2}}}{cx} \left[1 + \sum_{n=1}^N\frac{(2n-1)!!}{c^{2n}x^{2n}} + \overbrace{(2N+1)!! cx\frac{\int_{\nicefrac{cx}{\sqrt{2}}}^{cx}\frac{  \exp\paren*{\nicefrac{u^2}{2}}}{u^{2(N+1)}}du} {\exp\paren*{\nicefrac{c^2x^2}{2}}}}^A \right. \\ 
        &\left. -\overbrace{\sum_{n=0}^N\frac{(2n-1)!!\exp\paren*{-\nicefrac{c^2x^2}{4}}}{\sqrt{2^{-(2n+1)}}c^{2n}x^{2n}}}^{B} + \overbrace{cx\frac{\int_0^{\nicefrac{cx}{\sqrt{2}}} \exp\paren*{\nicefrac{u^2}{2}} du}{\exp\paren*{\nicefrac{c^2x^2}{2}}}}^C +  \right] 
    \end{align*}
    
    Notice that term (A) is of order $\OO{c^{-2(N+1)}}$ since
     \begin{align*}
       \lim_{c\rightarrow \infty} \frac{A}{c^{-2(N+1)}}
       & =  \lim_{c\rightarrow \infty}  (2N+1)!! x\frac{\int_{\nicefrac{cx}{\sqrt{2}}}^{cx} \frac{\exp\paren*{\nicefrac{u^2}{2}}}{u^{2(N+1)}} du} {\frac{\exp\paren*{\nicefrac{c^2x^2}{2}}}{c^{2N+3}}}\\
    &=\lim_{c\rightarrow \infty} (2N+1)!!x\frac{\frac{d}{dc}\int_{\nicefrac{cx}{\sqrt{2}}}^{cx} \frac{\exp\paren*{\nicefrac{u^2}{2}}}{u^{2(N+1)}} du}{\frac{d}{dc}\frac{\exp\paren*{\nicefrac{c^2x^2}{2}}}{c^{2N+3}}} \tag{applying \lopi}\\
         &=\lim_{c\rightarrow \infty} (2N+1)!!\frac{\frac{\exp\paren*{\nicefrac{c^2x^2}{2}}}{c^{2(N+1)}x^{2(N+1)}}\cdot x- \frac{\exp\paren*{\nicefrac{c^2x^2}{4}}2^{N+1}}{c^{2(N+1)}x^{2(N+1)}}\cdot \nicefrac{x}{\sqrt{2}}}{\frac{\exp\paren*{\nicefrac{c^2x^2}{2}} \paren*{(cx^2)c^{2N+3} - (2N+3)c^{2(N+1)}} }{c^{4N+6}}} \tag{applying Leibniz integral rule}\\
           &=\lim_{c\rightarrow \infty} (2N+1)!!\frac{c^{-2(N+1)}x^{-(2N+1)}\paren*{1- \exp\paren*{-\nicefrac{c^2x^2}{4}}\sqrt{2^{2N+1}}}}{c^{-2(N+1)}\paren*{x^2 - (2N+3)c^{-2}}}\\
          &=\lim_{c\rightarrow \infty} (2N+1)!!\frac{1}{x^{2N+1}}\frac{1- \exp\paren*{-\nicefrac{c^2x^2}{4}}\sqrt{2^{2N+1}}}{x^2 - (2N+3) c^{-2}}\\
         &= (2N+1)!!\frac{1}{x^{2N+3}},
    \end{align*}
    term (B) is $\OO{\exp\paren*{-\nicefrac{c^2x^2}{4}}}$ and so is term (C) since
       \begin{align*}
       \lim_{c\rightarrow \infty} \frac{C}{\exp\paren*{-\nicefrac{c^2x^2}{4}}}
       & =  \lim_{c\rightarrow \infty}  x\frac{\int_0^{\nicefrac{cx}{\sqrt{2}}} \exp\paren*{\nicefrac{u^2}{2}} du}{\frac{\exp\paren*{\nicefrac{c^2x^2}{4}}}{c}}\\
    &=\lim_{c\rightarrow \infty}x\frac{\frac{d}{dc}\int_0^{\nicefrac{cx}{\sqrt{2}}} \exp\paren*{\nicefrac{u^2}{2}} du}{\frac{d}{dc}\frac{\exp\paren*{\nicefrac{c^2x^2}{4}}}{c}} \tag{applying \lopi}\\
      &=\lim_{c\rightarrow \infty}x\frac{\exp\paren*{\nicefrac{c^2x^2}{4}}\cdot \frac{x}{\sqrt{2}}}{\frac{\exp\paren*{\nicefrac{c^2x^2}{4}} \paren*{(\nicefrac{cx^2}{2})c - 1} }{c^2}} \tag{applying Leibniz integral rule} \\
         &=\lim_{c\rightarrow \infty} \frac{\frac{x^2}{\sqrt{2}}}{\paren*{\nicefrac{x^2}{2} - \nicefrac{1}{c^2}}}\\
         &= \sqrt{2}.
    \end{align*}

   Consequently,
     \begin{align*}
        \frac{c\Im(cx)}{\exp\paren*{\nicefrac{c^2x^2}{2}}}  = \sqrt{\frac{2}{\pi}} \frac{1}{x}\brac*{ 1 +  \sum_{n=1}^N\frac{(2n-1)!!}{c^{2n}x^{2n}} + \OO{c^{-2(N+1)}} + \OO{\exp\paren*{-\nicefrac{c^2x^2}{4}}}} , 
    \end{align*}
    which proves Statement (\ref{item:Impoly}) and consequently Statement (\ref{item:Imlim}).
        Statement (\ref{item:Imlim}) implies that $\Im(cx)$ is $\OO{\frac{\exp\paren*{\nicefrac{c^2x^2}{2}}}{c}}$ when $x\neq 0$. Thus 
        $1+ \paren*{\Im(cx)}^2$ is $\OO{1}+\OO{\frac{\exp\paren*{c^2x^2}}{c^2}}$ and consequently $\OO{\exp\paren*{c^2x^2}}$ when $x\neq 0$. Notice that the $1+ \paren*{\Im(cx)}^2$ is trivially $\OO{\exp\paren*{c^2x^2}}$ when $x=0$ as well, which completes the proof of Statement (\ref{item:ImOO}).\\
        Statement (\ref{item:Imlim}) also implies that $\Im(cx)$ is $\MM{\frac{\exp\paren*{\nicefrac{c^2x^2}{2}}}{c}}$, when $x\neq 0$. Thus 
        $1+ \paren*{\Im(cx)}^2$ is $\MM{1}+\MM{\frac{\exp\paren*{c^2x^2}}{c^2}}$ and consequently $\MM{\frac{\exp\paren*{c^2x^2}}{c^2}}$, when $x\neq 0$. Notice that  $1+ \paren*{\Im(cx)}^2$ is trivially $\MM{\frac{\exp\paren*{c^2x^2}}{c^2}}$ when $x=0$ as well, which completes the proof of Statement (\ref{item:ImMM}).
\end{proof}
\end{lemma}

%\begin{theorem}
%    Let $A \neq 0$ be a $\dimm \times \dimm$ symmetric positive semi-definite matrix and $B$ be any matrix of the same dimension. Then there exists orthogonal column vectors $t_1$ and $t_2$ in  $\real^\dimm$ such that  $t'_1At_1>0$ and $B't_1$ and $B't_2$ are both standard basis vectors\textemdash with $1$ at one and only one position and $0$ at others\textemdash, distinct from each other.   
    
%\end{theorem}

%$\Lambda$ matrix for \uMSN. 
%$$\Lambda=\Omega^{\frac{1}{2}}\bar{\Delta}$$
%$$\bar{\Delta}= \bar{\Lambda}(I + \bar{\Lambda}'\bar{\Lambda})^{-\frac{1}{2}}$$
%As a consequence,
%$$I-\bar{\Delta}'\bar{\Delta}=(I+\bar{\Lambda}'\bar{\Lambda})^{-1}$$
%$$I-\bar{\Delta}\bar{\Delta}'=(I+\bar{\Lambda}\bar{\Lambda}')^{-1}$$

\section{Conclusions}
\label{sec:conclusions}
\noindent We give meaningful sufficient conditions that ensure identifiability of two\-/component mixtures with \SN, \rMSN, and \uMSN\ components. We proved identifiability in terms of the $\Gamma$ parameter that contains both the scale and the skewness information and has a consistent interpretation across the three skew normal families. Our results are strong in the sense that the set of parameter values not covered by the sufficient condition is a Lebesgue measure $0$ set in the parameter space. \cite{Ghosal2011} study the identifiability of a two-component mixture with the standard normal as one of the components and the second component itself given by a uncountable mixture of skew normals. Treating $G$ from their work as a point distribution, we can make a valid comparison between our identifiability result and theirs, concluding the superiority of our results, owing to a larger coverage of the parameter space by our conditions.

%
%The most relevant work \citet{ghosal2011identifiability} gives a similar theorem, however our theorem is more powerful.
%\include{algo}
\section*{\refname}
\bibliography{skewNormal,refs}{}

\begin{thebibliography}{29}
\providecommand{\natexlab}[1]{#1}
\providecommand{\url}[1]{\texttt{#1}}
\expandafter\ifx\csname urlstyle\endcsname\relax
  \providecommand{\doi}[1]{doi: #1}\else
  \providecommand{\doi}{doi: \begingroup \urlstyle{rm}\Url}\fi

\bibitem[Arellano-Valle and Genton(2005)]{arellano2005fundamental}
R.~B. Arellano-Valle and M.~G. Genton.
\newblock On fundamental skew distributions.
\newblock \emph{J Multivar Anal}, 96\penalty0 (1):\penalty0 93--116, 2005.

\bibitem[Azzalini(1985)]{azzalini1985class}
A.~Azzalini.
\newblock A class of distributions which includes the normal ones.
\newblock \emph{Scand J Stat}, 12\penalty0 (2):\penalty0 171--178, 1985.

\bibitem[Azzalini(1986)]{azzalini1986}
A.~Azzalini.
\newblock Further results on a class of distributions which includes the normal
  ones.
\newblock \emph{Statistica}, 46\penalty0 (2):\penalty0 199--208, 1986.

\bibitem[Azzalini and Dalla~Valle(1996)]{azzalini1996multivariate}
A.~Azzalini and A.~Dalla~Valle.
\newblock The multivariate skew-normal distribution.
\newblock \emph{Biometrika}, 83\penalty0 (4):\penalty0 715--726, 1996.

\bibitem[Blanchard et~al.(2010)Blanchard, Lee, and Scott]{Blanchard2010}
G.~Blanchard, G.~Lee, and C.~Scott.
\newblock Semi-supervised novelty detection.
\newblock \emph{J Mach Learn Res}, 11:\penalty0 2973--3009, 2010.

\bibitem[Bordes et~al.(2006)Bordes, Delmas, and
  Vandekerkhove]{bordes2006semiparametric}
L.~Bordes, C.~Delmas, and P.~Vandekerkhove.
\newblock Semiparametric estimation of a two-component mixture model where one
  component is known.
\newblock \emph{Scand J Stat}, 33\penalty0 (4):\penalty0 733--752, 2006.

\bibitem[Browne and McNicholas(2015)]{browne2015mixture}
R.~P. Browne and P.~D. McNicholas.
\newblock A mixture of generalized hyperbolic distributions.
\newblock \emph{Can J Stat}, 43\penalty0 (2):\penalty0 176--198, 2015.

\bibitem[Dempster et~al.(1977)Dempster, Laird, and Rubin]{dempster1977maximum}
A.~P. Dempster, N.~M. Laird, and D.~B. Rubin.
\newblock Maximum likelihood from incomplete data via the {EM} algorithm.
\newblock \emph{J R Stat Soc B}, pages 1--38, 1977.

\bibitem[Genton(2004)]{genton2004skew}
M.~G. Genton.
\newblock \emph{Skew-elliptical distributions and their applications: a journey
  beyond normality}.
\newblock CRC Press, 2004.

\bibitem[Ghosal and Roy(2011)]{Ghosal2011}
S.~Ghosal and A.~Roy.
\newblock Identifiability of the proportion of null hypotheses in skew-mixture
  models for the p-value distribution.
\newblock \emph{Electron J Statist}, 5:\penalty0 329--341, 2011.

\bibitem[Holzmann et~al.(2006)Holzmann, Munk, and
  Gneiting]{holzmann2006identifiability}
H.~Holzmann, A.~Munk, and T.~Gneiting.
\newblock Identifiability of finite mixtures of elliptical distributions.
\newblock \emph{Scand J Stat}, 33\penalty0 (4):\penalty0 753--763, 2006.

\bibitem[Jain et~al.(2016{\natexlab{a}})Jain, White, and Radivojac]{Jain2016b}
S.~Jain, M.~White, and P.~Radivojac.
\newblock Estimating the class prior and posterior from noisy positives and
  unlabeled data.
\newblock In \emph{Advances in Neural Information Processing Systems}, NIPS
  2016, pages 2693--2701, 2016{\natexlab{a}}.

\bibitem[Jain et~al.(2016{\natexlab{b}})Jain, White, Trosset, and
  Radivojac]{Jain2016}
S.~Jain, M.~White, M.~W. Trosset, and P.~Radivojac.
\newblock Nonparametric semi-supervised learning of class proportions.
\newblock \emph{arXiv preprint arXiv:1601.01944}, 2016{\natexlab{b}}.
\newblock URL \url{http://arxiv.org/abs/1601.01944}.

\bibitem[Jain et~al.(2017)Jain, White, and Radivojac]{Jain2017}
S.~Jain, M.~White, and P.~Radivojac.
\newblock Recovering true classifier performance in positive-unlabeled
  learning.
\newblock In \emph{Proceedings of the 31st AAAI Conference on Artificial
  Intelligence}, AAAI 2017, pages 2066--2072, 2017.

\bibitem[Kim and Genton(2011)]{kim2011characteristic}
H.~M. Kim and M.~G. Genton.
\newblock Characteristic functions of scale mixtures of multivariate
  skew-normal distributions.
\newblock \emph{J Multivar Anal}, 102\penalty0 (7):\penalty0 1105--1117, 2011.

\bibitem[Lee and McLachlan(2013)]{lee2013mixtures}
S.~X. Lee and G.~J. McLachlan.
\newblock On mixtures of skew normal and skew t-distributions.
\newblock \emph{Adv Data Anal Classif}, 7\penalty0 (3):\penalty0 241--266,
  2013.

\bibitem[Lin(2009)]{lin2009maximum}
T.~I. Lin.
\newblock Maximum likelihood estimation for multivariate skew normal mixture
  models.
\newblock \emph{J Multivar Anal}, 100\penalty0 (2):\penalty0 257--265, 2009.

\bibitem[McLachlan and Peel(2000)]{McLachlan2000}
G.~J. McLachlan and D.~Peel.
\newblock \emph{Finite mixture models}.
\newblock John Wiley \& Sons, Inc., 2000.

\bibitem[Menon et~al.(2015)Menon, van Rooyen, Ong, and Williamson]{Menon2015}
A.~K. Menon, B.~van Rooyen, C.~S. Ong, and R.~C. Williamson.
\newblock Learning from corrupted binary labels via class-probability
  estimation.
\newblock In \emph{Proceedings of the 32nd International Conference on Machine
  Learning}, ICML 2015, pages 125--134, 2015.

\bibitem[Patra and Sen(2016)]{Patra2016}
R.~K. Patra and B.~Sen.
\newblock Estimation of a two-component mixture model with applications to
  multiple testing.
\newblock \emph{J R Statist Soc B}, 78\penalty0 (4):\penalty0 869--893, 2016.

\bibitem[Pewsey(2003)]{pewsey2003characteristic}
A.~Pewsey.
\newblock The characteristic functions of the skew-normal and wrapped
  skew-normal distributions.
\newblock In \emph{27 Congreso Nacional de Estadistica e Investigaci{\'o}n
  Operativa}, pages 4383--4386, 2003.

\bibitem[Storey(2002)]{Storey2002}
J.~D. Storey.
\newblock A direct approach to false discovery rates.
\newblock \emph{J R Stat Soc B}, 64\penalty0 (3):\penalty0 479--498, 2002.

\bibitem[Storey(2003)]{Storey2003a}
J.~D. Storey.
\newblock The positive false discovery rate: a {Bayesian} interpretation and
  the q-value.
\newblock \emph{Ann Stat}, 31\penalty0 (6):\penalty0 2013--2035, 2003.

\bibitem[Storey and Tibshirani(2003)]{Storey2003}
J.~D. Storey and R.~Tibshirani.
\newblock Statistical significance for genomewide studies.
\newblock \emph{Proc Natl Acad Sci U S A}, 100\penalty0 (16):\penalty0
  9440--9445, 2003.

\bibitem[Tallis(1961)]{tallis1961moment}
G.~M. Tallis.
\newblock The moment generating function of the truncated multi-normal
  distribution.
\newblock \emph{J R Stat Soc B}, pages 223--229, 1961.

\bibitem[Tallis and Chesson(1982)]{Tallis1982}
G.~M. Tallis and P.~Chesson.
\newblock Identifiability of mixtures.
\newblock \emph{J Austral Math Soc Ser A}, 32:\penalty0 339--348, 1982.

\bibitem[Teicher(1963)]{teicher1963identifiability}
H.~Teicher.
\newblock Identifiability of finite mixtures.
\newblock \emph{Ann Math Stat}, 34\penalty0 (4):\penalty0 1265--1269, 1963.

\bibitem[Ward et~al.(2009)Ward, Hastie, Barry, Elith, and Leathwick]{Ward2009}
G.~Ward, T.~Hastie, S.~Barry, J.~Elith, and J.~R. Leathwick.
\newblock Presence-only data and the {EM} algorithm.
\newblock \emph{Biometrics}, 65\penalty0 (2):\penalty0 554--563, 2009.

\bibitem[Yakowitz and Spragins(1968)]{Yakowitz1968}
S.~J. Yakowitz and J.~D. Spragins.
\newblock On the identifiability of finite mixtures.
\newblock \emph{Ann Math Statist}, 39:\penalty0 209--214, 1968.

\end{thebibliography}
\bibliographystyle{plainnat}
\newpage
\appendix
\section{\uMSN\ characteristic function}
\label{sec:CFUSNCF}
\begin{theorem}
\label{thm:CFUSNCF}
The characteristic function of \uMSN$(\mu,\Omega,\Lambda)$ is given by
$$\cf(t)=\exp\braces*{\img t'\mu - \frac{1}{2}t'\Omega t}\prod_{i=1}^\dimm\paren*{1+\img\Im(\Lambda_i't)},$$
where $\dimm$ is the dimensionality of the \uMSN\ distribution, $\img$ is the imaginary number, $\Im(x)=\int_0^x \sqrt{\nicefrac{2}{\pi}} \exp\paren*{\nicefrac{u^2}{2}}du$ and $\Lambda_i$ is the $i$\textsuperscript{th} column of $\Lambda$. 
\begin{proof}
We use the stochastic representation of $X \sim CFUSN\paren*{\mu,\Omega,\Lambda}$ obtained from \cite{lin2009maximum}, given by $X=\Lambda H + G$, where $H\sim \TN\paren*{0,\iden,\orthant}$, the standard multivariate normal distribution truncated below $0$ in all the dimensions and $G \sim N_\dimm(\mu,\Gamma)$, for $\Gamma=\Omega-\Lambda\Lambda'$ \textemdash a symmetric positive semi-definite matrix. It follows that the \cff\  of $X$ can be expressed in terms of \cff's of Normal distribution and truncated Normal distribution; precisely, $\cf_X(t)=\cf_{G}(t) \cdot \cf_{\Lambda H}(t)$. Using the expression for the \cf\ of Multivariate Normal,
\begin{align}
    \cf_X(t)=\exp\braces*{\img t'\mu - \frac{1}{2}t'\Gamma t}  \cf_{\Lambda H}(t) \label{eq:cfProd}
\end{align}
Basic properties of a \cff\ and its connection with the corresponding \mgf\ gives $\cf_{\Lambda H}(t)=\cf_H(\Lambda't)=\mgf_H(\img \Lambda't)$.  Using the expression for $\mgf_H(t)$, derived in \cite{tallis1961moment} (p. 225), and replacing $t$ by $\img\Lambda' t$  and $R$ (the covariance matrix in \cite{tallis1961moment})  by $\iden$ ($\dimm \times \dimm$ identity matrix), we get
\begin{align*}
CF_{\Lambda H}(t)&=\exp\braces*{-\frac{1}{2}t'\Lambda\Lambda't}\frac{\paren*{2\pi}^{-\nicefrac{\dimm}{2}}\int_{\orthant} \exp\braces*{-\frac{1}{2}(w-\img\Lambda' t)'(w-\img\Lambda' t)}dw}{\int_{\orthant} \phi_\dimm(u)du}  \tag{where $\phi_\dimm$ is the pdf of \normal$(0,\iden)$}\\
&=  2^\dimm\exp\braces*{-\frac{1}{2}t'\Lambda\Lambda't} \int_{\orthant}\paren*{2\pi}^{-\nicefrac{\dimm}{2}}\exp\braces*{-\frac{1}{2}\sum_{i=1}^\dimm(w_i-\img  \Lambda_i't)^2}dw \tag{where $w=[w_i]_{i=1}^\dimm$}\\
&=  2^\dimm\exp\braces*{-\frac{1}{2}t'\Lambda\Lambda't}
\int_{\orthant}\paren*{2\pi}^{-\nicefrac{\dimm}{2}}\prod_{i=1}^\dimm\exp\braces*{-\frac{1}{2}(w_i-\img \Lambda_i't)^2}dw\\
&= 2^\dimm\exp\braces*{-\frac{1}{2}t'\Lambda\Lambda't}\prod_{i=1}^\dimm\int_{0}^{\infty} \paren*{2\pi}^{-\nicefrac{1}{2}}\exp\braces*{-\frac{1}{2}(w_i-\img\Lambda_i't)^2}dw_i.\\
\end{align*}
Applying the substitution $u_i=-w_i+\img \Lambda_i't$ for the integral in the numerator, changes the domain of the integration from  the real line to the complex plane. To define such an integral correctly, one needs to specify the path in the complex plane across which the integration is performed. Using the path from $-\infty +\img \Lambda_i't$ to $\img \Lambda_i't$, parallel to the real line, we get  
\begin{align*}
CF_{\Lambda H}(t)&=2^\dimm\exp\braces*{-\frac{1}{2}t'\Lambda\Lambda't} \prod_{i=1}^\dimm\int_{-\infty+\img \Lambda_i't}^{\img \Lambda_i't} \paren*{2\pi}^{-\nicefrac{1}{2}} \exp\braces*{-\frac{1}{2}u_i^2}du_i \\
&=2^\dimm\exp\braces*{-\frac{1}{2}t'\Lambda\Lambda't} \prod_{i=1}^\dimm\int_{-\infty+\img \Lambda_i't}^{\img \Lambda_i't} \phi(u_i)du_i.
\end{align*}
Using Lemma 1 from \cite{kim2011characteristic} to simplify the integral term, we get
\begin{align*}
CF_{\Lambda H}(t) &= 2^\dimm \exp\braces*{-\frac{1}{2}t'\Lambda\Lambda't} \prod_{i=1}^\dimm \paren*{\frac{1}{2} + \img \frac{1}{\sqrt{\pi}} \int_0^{\frac{\Lambda_i't}{\sqrt{2}}} \exp\braces*{u_i^2}du_i}\\
&= 2^\dimm\exp\braces*{-\frac{1}{2}t'\Lambda\Lambda't} \prod_{i=1}^\dimm \paren*{\frac{1}{2} + \img \frac{1}{2}  \int_0^{\Lambda_i't} \sqrt{\frac{2}{\pi}} \exp\braces*{\frac{v_i^2}{2}}dv_i}\tag{substituting $v_i=\nicefrac{u_i}{\sqrt{2}}$}\\
&= \exp\braces*{-\frac{1}{2}t'\Lambda\Lambda't} \prod_{i=1}^\dimm \paren*{1 + \img \Im\paren*{\Lambda_i't}},
\end{align*}
Substituting the expression for $CF_{\Lambda H}(t)$ in equation (\ref{eq:cfProd}) completes the proof.

%\begin{align*}
%CF^{HN}_{\Lambda}(ct)&= \frac{\int_{-ic\Lambda't}^{\infty}\prod_{i=1}^k\exp\braces*{-\frac{1}{2} u_i^2}du_1\ldots du_k}{\int_{0}^{\infty}\prod_{i=1}^k\exp\braces*{-\frac{1}{2} u_i^2}du_1\ldots du_k}\exp\braces*{-\frac{1}{2}c^2t'\Lambda\Lambda't} \\
%&=\frac{\int_{-ic\Lambda't}^{\infty}\exp\braces*{-\frac{1}{2}u'u}du}{\int_{0}^{\infty}\exp\braces*{-\frac{1}{2}u'u}du}\exp\braces*{-\frac{1}{2}c^2t'\Lambda\Lambda't} \\
%\end{align*}
%The characteristic function of the normal term, $U$, is given by
%$$CF^{N}(t)=\exp\braces*{\img\mu't - \frac{1}{2}t'\Gamma t}.$$
%The characteristic function of the CFUSN distribution is given by the product of $CF^{N}(t)$ and $CF^{HN}(\Lambda't)$;
%\begin{align*}
%    CF(t)&=\exp\braces*{\img\mu't - \frac{1}{2}t'(\Gamma+\Lambda\Lambda') t}\prod_{i=1}^\dimm \paren*{1 + \img \Im\paren*{\Lambda_i't}}\\
%    &=\exp\braces*{\img\mu't - \frac{1}{2}t'\Omega t}\prod_{i=1}^\dimm \paren*{1 + i \Im\paren*{\Lambda_i't}},
%\end{align*}
%where $\Omega=\Gamma + \Lambda\Lambda'$
\end{proof}
\end{theorem}

\end{document}